\documentclass[11pt,twoside]{article}

\usepackage{epsfig,amsfonts,color}
\usepackage[final]{pdfpages}
\usepackage{amsmath}
\usepackage{blkarray}
\usepackage{bm}
\usepackage{algorithm}
\usepackage{algpseudocode}
\renewcommand{\algorithmicfunction}{\textbf{Function}}
\algdef{SE}[TASK]{Task}{EndTask}%
   [2]{\algorithmicfunction\ \tt{#1}\ifthenelse{\equal{#2}{}}{}{(#2)}}%
   {\algorithmicend\ \algorithmicfunction}%

\usepackage{graphicx}
\usepackage[margin = 1cm]{caption}
\usepackage{subcaption}
\usepackage{placeins}

\newcommand{\mods}[1]{\left|  #1 \right|}

\usepackage{amssymb, palatino, geometry,url}

\usepackage[colorlinks=true,linkcolor=blue,citecolor=blue,urlcolor=blue]{hyperref}
\usepackage[titletoc,title]{appendix}
\usepackage{multirow}
\usepackage{mathtools}
\usepackage{fancyhdr}
\newcommand{\be}{\begin{equation}}
\newcommand{\ee}{\end{equation}}
\newcommand{\bea}{\begin{eqnarray}}
\newcommand{\eea}{\end{eqnarray}}

\newcommand{\bvec}{\left(\begin{array}{c}}
\newcommand{\evec}{\end{array}\right)}
\newcommand{\bsub}{\begin{subequations}}
\newcommand{\esub}{\end{subequations}}

\usepackage{accsupp}

\usepackage{listings}
\DeclareCaptionFont{white}{\color{white}}
\DeclareCaptionFormat{listing}{\colorbox[cmyk]{0.7, 0.35, 0.35,0.01}{\parbox{\dimexpr\textwidth-2\fboxsep\relax}{#1#2#3}}}
\usepackage[dvipsnames]{xcolor}

\lstdefinelanguage{julia}
{
keywordsprefix=\@,
morekeywords={
exit,whos,edit,load,is,isa,isequal,typeof,tuple,ntuple,uid,hash,finalizer,convert,promote,
subtype,typemin,typemax,realmin,realmax,sizeof,eps,promote_type,method_exists,applicable,
invoke,dlopen,dlsym,system,error,throw,assert,new,Inf,Nan,pi,im,begin,while,for,in,return,
break,continue,macro,quote,let,if,elseif,else,try,catch,end,bitstype,ccall,do,using,module,
import,export,importall,baremodule,immutable,local,global,const,Bool,Int,Int8,Int16,Int32,
Int64,Uint,Uint8,Uint16,Uint32,Uint64,Float32,Float64,Complex64,Complex128,Any,Nothing,None,
function,type,typealias,abstract,get_node,add_edge,create_estimation_model,set_solution,
solve, get_solution, solve_ss_problem, create_estimation_problem, addnode
},
morekeywords = [2]{triggered_by,compute_time,trigger_during_busy,send_on,delay,send_wait,start},
sensitive=true,
morecomment=[l]{\#},
morestring=[b]',
morestring=[b]"
}
\usepackage{authblk}
\usepackage{amsthm}
\theoremstyle{plain}

\begin{document}

\title{Graph-Based Modeling and Simulation of Complex Systems}

 \author{Jordan Jalving$^{\dag \ddag}$, Yankai Cao$^{\dag}$ and Victor M. Zavala$^{\dag}$\\
  {\small $^{\dag}$ Department of Chemical and Biological Engineering}\\
  {\small University of Wisconsin-Madison, 1415 Engineering Dr, Madison, WI 53706, USA}\\
  {\small $^{\ddag}$ Decision and Infrastructure Sciences Division}\\
  {\small Argonne National Laboratory, 9700 South Cass Ave, Lemont, IL 60439, USA}\\
 }

\date{}

\maketitle

\abstract {We present graph-based modeling abstractions to represent cyber-physical dependencies arising in complex systems.
Specifically, we propose an algebraic graph abstraction to capture physical connectivity in complex optimization models and a computing graph
abstraction to capture communication connectivity in computing architectures.
The proposed abstractions are scalable and are used as the backbone of a ${\tt Julia}$-based software package that we call ${\tt Plasmo.jl}$.
We show how the {algebraic graph} abstraction facilitates the implementation, analysis, and decomposition of optimization problems and we show how the {computing graph}
abstraction facilitates the implementation of optimization and control algorithms and their simulation in virtual environments that involve distributed, centralized, and hierarchical computing architectures.
}
\\

{\bf Keywords}: graphs, cyber-physical, connectivity, algebraic, computing

\section{Introduction}

Complex systems are {\em cyber-physical} in nature, in the sense that a physical system is driven by decisions made by a cyber (computing) system \cite{lee2015cyber}.
For instance, a chemical process is a physical system that is driven by decisions made by a control system, which is in turn a cyber system comprised of a collection of devices (e.g., sensors, controllers, actuators)
that execute diverse computing tasks (e.g., data processing and control action computation) and that exchange signals and data (e.g., measurements and actions) through a communication network.
The devices executing the tasks of a control system form a {\em computing architecture}, similar in spirit to a parallel computing cluster in which computing processors are connected through a communication network.
Other examples of cyber-physical systems include hierarchical architectures for coordination of supply chain, scheduling, and planning tasks for an enterprise or the control architecture for an infrastructure network such as a natural gas pipeline \cite{Brunaud2017,baldea_daoutidis_2012}.

{Modeling and simulating the behavior of cyber-physical systems is becoming increasingly important but also increasingly challenging}.
In particular, emerging paradigms such as cloud computing and the internet-of-things are drastically changing the landscape of decision-making architectures.
Such reconfiguration is driven by the need to process increasingly larger amounts of data in a distributed manner while making decisions faster and in a more scalable manner.
Architectural reconfiguration needs to carefully balance diverse cyber-physical issues such as economic performance, safety, data privacy, as well as computing and communication latency and failures.
For instance, failure of a computing device (e.g., a sensor) can lead to significant losses in performance or to full collapse of the physical system.

Capturing interdependencies between cyber and physical systems in a coherent manner is technically challenging.  Specifically, the behavior of a physical system is expressed mathematically in the form of an {\em algebraic model} (i.e., a set of algebraic equations) while the behavior of a cyber system is expressed mathematically in the form of {\em algorithms}.  Moreover, in modeling a cyber system, one must consider the fact that algorithms are executed under highly heterogeneous and dynamic computing architectures that exhibit complex computing and communication protocols and logic (e.g., synchronous and asynchronous) and associated time delays \cite{lee2008cyber}.

This work proposes \emph{graph-based} abstractions to facilitate modeling and simulation of cyber-physical systems.
The proposed graph-based abstractions seek to provide a coherent framework to capture modeling elements that arise in various engineering domains.
In the context of modeling physical systems (specifically chemical processes), early sequential modular and equation-based flowsheeting tools used graph-theoretical insights to express equations in a modular form and to
facilitate object-oriented software implementations, analysis, and algorithmic development \cite{Westerberg1994}.
Such developments are the basis of powerful simulation environments such as ${\tt Ascend}$, ${\tt AspenPlus}$, and ${\tt gPROMS}$ \cite{gProms,Piela1990}.
Graph-theoretical concepts have also been widely used for expressing and processing algebraic models in platforms such as
${\tt AMPL}$, ${\tt GAMS}$, ${\tt AIMMS}$, and ${\tt CasADi}$ \cite{Fourer1990,GamsSoftware2013,AIMMS,Andersson2018}.
The Modelica \cite{Fritzson02modelica} simulation platform also uses graph concepts such as modularity and inheritance to
instantiate and simulate complex systems.
Recent advances in algebraic modeling languages have been enabled by ${\tt Pyomo}$ \cite{hart2017pyomo} and ${\tt JuMP}$ \cite{DunningHuchetteLubin2017}.
These open-source tools enable the user to implement models in modular form and to expose structure to facilitate the implementation of parallel solution algorithms.
Graph-based modeling abstractions have recently been explored in convex optimization \cite{Hallac2015}, infrastructure networks \cite{Jalving2017}, and simulation of partial differential equations  \cite{abhyankar2018petsc},
but these abstractions are restricted in that the graph structure is directly tied to physical topology, thus limiting modeling flexibility.

In the context of modeling and simulating cyber systems, the most popular framework is ${\tt Simulink}$ \cite{Simulink}. Here, a system is expressed in terms of blocks of operators (computing tasks) which are
connected by communication channels and input and output ports.  This framework facilitates the simulation of complex control architectures.
Agent-based modeling platforms can also be used to simulate cyber systems; under this abstraction, agents make decisions and communicate under channels \cite{Macal2005}.
Popular agent-based simulation tools include ${\tt RePast}$, ${\tt MASON}$, and ${\tt Swarm} $ \cite{North2013,Collier2013,LukeMASON,minar1996swarm}.

While the aforementioned physical and cyber modeling tools exploit graph concepts to facilitate implementation, they do not use a coherent abstraction (which is key to enable extensibility).
To enable this, in this work, we introduce general abstractions for modeling cyber-physical systems.  Specifically, we propose the concept of an \emph{algebraic graph} to facilitate modeling of mathematical optimization problems and
the concept of a \emph{computing graph} to facilitate simulation of cyber systems. The graph abstractions exploit
\emph{physical and communication topologies} to facilitate model construction, data management, and analysis.
The computing graph abstraction offers advantages over ${\tt Simulink}$ and agent-based tools in that it handles cyber features such as communication delays, latency, and synchronous/asynchronous computing
and information exchange in a more coherent manner. This is done by using a {\em state-space} representation of the computing graph, which keeps track of task states and which manages
actions and timings that trigger state changes in the form of action signals.  We also demonstrate how the proposed abstractions can be combined to enable modeling of cyber-physical systems and we
provide an implementation in the ${\tt Julia}$ programming language \cite{Bezanson2017} that we call ${\tt Plasmo.jl}$.

The manuscript is structured as follows.  Section \ref{sec:graph_representations} presents basic graph terminology for describing
the proposed abstractions. Section \ref{sec:graph_modeling} presents the algebraic graph abstraction for modeling optimization problems
with a focus on facilitating data management and problem decomposition.  Section \ref{sec:graph_modeling} also introduces
our own implementation of graph abstractions called ${\tt Plasmo.jl}$ with a simple example.
Section \ref{sec:workflows} introduces the computing graph abstraction and
Section \ref{sec:case_study} provides case studies for cyber-physical
systems.

\section{Basic Graph Terminology}\label{sec:graph_representations}

\begin{figure}[ht]
    \centering
    \includegraphics[scale=0.4]{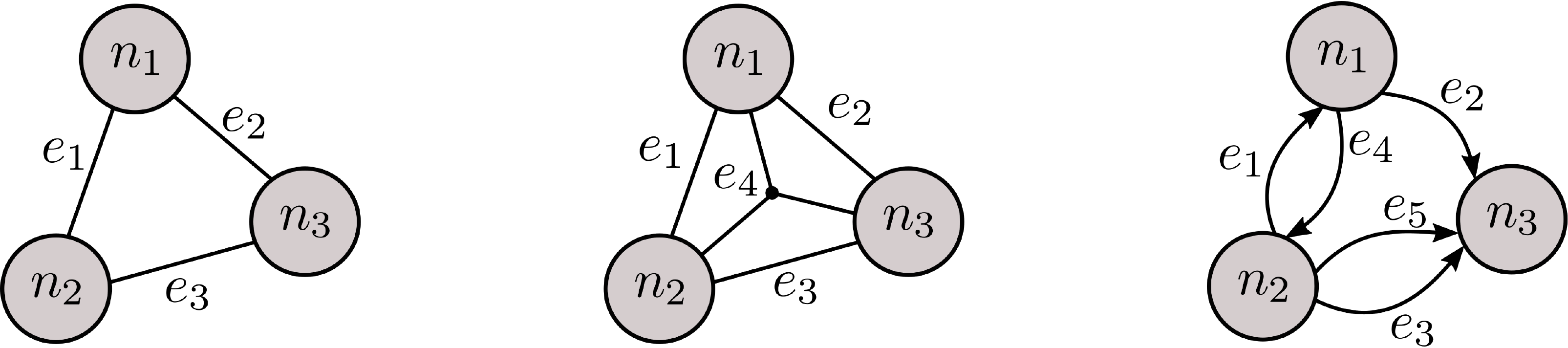}
    \caption{Example graph representations.  A simple undirected graph (left) with three nodes and three simple edges. An undirected hypergraph (center)
    with three nodes, three edges and one hyperedge.  A directed multigraph (right) with three nodes and five directed edges.
    }
    \label{fig:graph_representations}
\end{figure}

A graph $\mathcal{G}(\mathcal{N},\mathcal{E})$ is a collection of nodes $\mathcal{N}$ and edges $\mathcal{E}$.
We highlight the fact that a set of nodes belong to a specific graph $\mathcal{G}$ by using the syntax $\mathcal{N}(\mathcal{G})$ and we denote the node elements using 
index $n\in\mathcal{N}(\mathcal{G})$. Similarly, we highlight that a set of edges belong to $\mathcal{G}$ by using syntax
$\mathcal{E}(\mathcal{G})$ and we denote the edge elements using $e\in \mathcal{E}(\mathcal{G})$.  We define the set of supporting nodes of an edge $e$ (nodes that the edge connects)
as $\mathcal{N}(e)\subseteq\mathcal{N}$ and the
set of supporting edges for a node $n$ (edges connected to the node) as $\mathcal{E}(n)\subseteq \mathcal{E}$.
In a standard graph, two nodes are connected by an edge $e$.  For example, nodes $n_1$ and $n_2$ are connected by the single
edge $e_1$ in the left panel of Figure \ref{fig:graph_representations}.  In a hypergraph, multiple nodes can be connected by a single edge.  Hypergraphs are useful for
describing algebraic structures of systems. For example, the middle panel in Figure \ref{fig:graph_representations} is a hypergraph that contains an edge $e_4$ that connects all three nodes. In a multigraph, multiple edges can connect two nodes ($\mathcal{E}(n)$ might not be a singleton). In a directed multigraph, multiple edges can connect two nodes and also have direction.  This representation is useful when edges represent flows (e.g., physical or communication) between nodes.  The set of incoming directed edges to a node $n$ is expressed as $\mathcal{E}_{in}(n)$ and the set of outgoing edges is $\mathcal{E}_{out}(n)$. The right-most graph in Figure \ref{fig:graph_representations} contains a directed multigraph wherein the nodes are connected by directed edges.

Connectivity between nodes and edges is usually expressed in terms of the incidence matrix $A\in \mathcal{R}^{|\mathcal{N}|\times |\mathcal{E}|}$ (where
the notation $|\mathcal{S}|$ for set $\mathcal{S}$ denotes its cardinality). For an undirected standard graph we have $A_{n,e}=1$ if $n\in\mathcal{N}(e)$ or $e\in\mathcal{E}(n)$.
For a directed graph we have $A_{n,e}=1$ if $e\in\mathcal{E}_{in}(n)$ and $A_{n,e}=-1$ if $e\in\mathcal{E}_{out}(n)$.
Moreover, for a standard directed graph we have that $\sum_{n\in\mathcal{N}}A_{n,e}=0$ for all $e\in\mathcal{E}$ and for a
standard undirected graph we have $\sum_{n\in\mathcal{N}}A_{n,e}=2$ for all $e\in\mathcal{E}$ (because in a standard graph an edge only has two supporting nodes).
For a standard graph, the degree of a node $n$ is the number of edges connected to the node and can be computed as $\textrm{deg}(n)=|\mathcal{E}(n)|=\sum_{e\in\mathcal{E}}|A_{n,e}|$. The set of nodes connected to node $n$ (without counting self-loops) is denoted as $\mathcal{N}(n)$ and for a standard graph we have that $|\mathcal{N}(n)|=|\mathcal{E}(n)|$.

A graph can be constructed in a {\em hierarchical manner} by partitioning it into subgraphs (with their own local nodes and edges). The nodes in a subgraph are connected to those in other subgraphs by a set of global edges (see Figure \ref{fig:subgraph_structures}).
Consider, for instance, the subgraphs $\mathcal{G}_i,\mathcal{G}_j$ and the graph
$\mathcal{G}(\{\mathcal{G}_i,\mathcal{G}_j\}, \mathcal{E})$,  which indicates that $\mathcal{G}=(\mathcal{N},\mathcal{E})$ is a graph
with nodes $\mathcal{N}(\mathcal{G}_i)\cup \mathcal{N}(\mathcal{G}_j)$ and edges $\mathcal{E}=\mathcal{E}_g\cup\mathcal{E}(\mathcal{G}_i)\cup \mathcal{E}(\mathcal{G}_j)$ and where,
for every global edge $e\in\mathcal{E}_g$, we have that $\mathcal{N}(e)\in \mathcal{E}(\mathcal{G}_i)$ and $\mathcal{E}(\mathcal{G}_j)$.
In other words, the edges in $\mathcal{E}_g$ only connect nodes across subgraphs $\mathcal{G}_i$ and $\mathcal{G}_j$ but not within subgraphs.
Consequently, if we treat the elements of a node set $\mathcal{N}$ as graphs, we can represent  $\mathcal{G}(\{\mathcal{G}_i,\mathcal{G}_j\}, \mathcal{E})$ as a general
graph of the form $\mathcal{G}(\mathcal{N}, \mathcal{E})$. This nesting procedure can be repeated over multiple levels to form a hierarchical (multi-level) graph.

The discussed terminology and constructs form the basis for a wide range of graph analysis techniques.  Graphs also provide a coherent framework to represent and analyze a wide range of complex systems. 

\section{Algebraic Graphs}\label{sec:graph_modeling}

Here we present a graph-based modeling abstraction to enable the flexible representation of algebraic optimization models.

\subsection{Representation}

We define an \textbf{algebraic model graph} $\mathcal{MG}(\mathcal{N},\mathcal{E})$ (that we simply refer to as a model graph) as a {\em hypergraph} wherein
every node $n\in\mathcal{N}(\mathcal{MG})$ has an associated algebraic {\em component} model of the form:
    \begin{align}\label{eq:subprobs}
    {\min}_{x_n\in\mathcal{X}_n} & \quad f_n(x_n),
    \end{align}
where $x_n$ is a vector of decision variables (of arbitrary dimension), $\mathcal{X}_{n}:=\{x\,|\,c_n(x) \geq 0 \}$ is the feasible set with associated constraint vector mapping $c_n(\cdot)$, and the mapping $f_n(\cdot)$ is a scalar objective function.  For simplicity, we define the elements of the node  $\mathcal{N}$ as $\{n_1,n_2,...,n_{|\mathcal{N}|}\}$.  We highlight that this general representation also includes weighted graphs, systems of algebraic equations ,and optimization problems with mixed-integer variables. 

To capture coupling between the component models residing in the nodes, we introduce the concept of \emph{link constraints}. For simplicity in the presentation, here we consider linear linking constraints of the form:
\begin{align}\label{eq:link_constraint}
\sum_{n \in \mathcal{N}(e)} \Pi_{e,n} x_n = 0 \quad e \in \mathcal{E}(\mathcal{MG}),
\end{align}
where $e \in \mathcal{E}(\mathcal{MG})$ are the hyperedges of the model graph and $\mathcal{N}(e)$ is the set of nodes that support edge $e$ (nodes connected to the edge). The connectivity matrix $\Pi_{e,n}$ corresponds to the coefficients of the linking constraint.   For simplicity, we define the elements of the edge set  $\mathcal{E}$ as $\{e_1,e_2,...,e_{|\mathcal{E}|}\}$.

To enable compact notation, we can encapsulate all the link constraints as $\Pi_{\mathcal{MG}}x_{\mathcal{MG}}=0$, where $x_{\mathcal{MG}}$ denotes the set of all decision variables in the model graph $\mathcal{MG}$ and where $\Pi_{\mathcal{MG}}$ is a connectivity matrix of the form:
\begin{align}
\Pi_{\mathcal{MG}} =\left(
\begin{array}{ccccc}
 \Pi_{e_1,n_1} &\Pi_{e_1,n_2} & \dotsm & \Pi_{e_1,n_{|\mathcal{N}|}}  \\
 \Pi_{e_2,n_1} &\Pi_{e_2,n_2} & \dotsm & \Pi_{e_2,n_{|\mathcal{N}|}}  \\
 \vdots & \vdots & \vdots & \vdots  \\
  \Pi_{e_{|\mathcal{E}|},n_1}  & \Pi_{e_{|\mathcal{E}|},n_2} & \dotsm & \Pi_{e_{|\mathcal{E}|},n_{|\mathcal{N}|}}  \\
\end{array}\right).
 \end{align}
Under this representation, we have that $\Pi_{e_i,n_j}=0$ if $n_j\notin \mathcal{N}(e_i)$. In other words, the block matrix of node $n_j$ is zero if the node does not support the edge $e_i$ (it is not connected to it). Consequently, the connectivity matrix $\Pi_{\mathcal{MG}}$ captures connectivity over the entire model graph. In fact, the node-edge sparsity structure coincides with that of the incidence matrix.

We use the above notation to express the optimization problem over the entire graph $\mathcal{MG}$ as:
\begin{subequations}\label{eq:model-graph-compact}
    \begin{align}
        \min_{x_\mathcal{MG}} & \quad \sum_{n \in \mathcal{N}(\mathcal{MG})} f_n(x_n) \\
        \textrm{s.t.} & \quad x_n \in \mathcal{X}_n,   \quad n \in \mathcal{N}(\mathcal{MG}) \\
        & \quad \Pi_{\mathcal{MG}} x_{\mathcal{MG}} = 0.   \label{eq:link-constraint}
    \end{align}
\end{subequations}
We have assumed, for convenience, that the objectives of every node are added to form the graph objective. Other operators can be performed on the aggregate the objectives (e.g., product, expected value, variance, worst-case).

\begin{figure}[ht]
    \centering
    \includegraphics[scale=0.17]{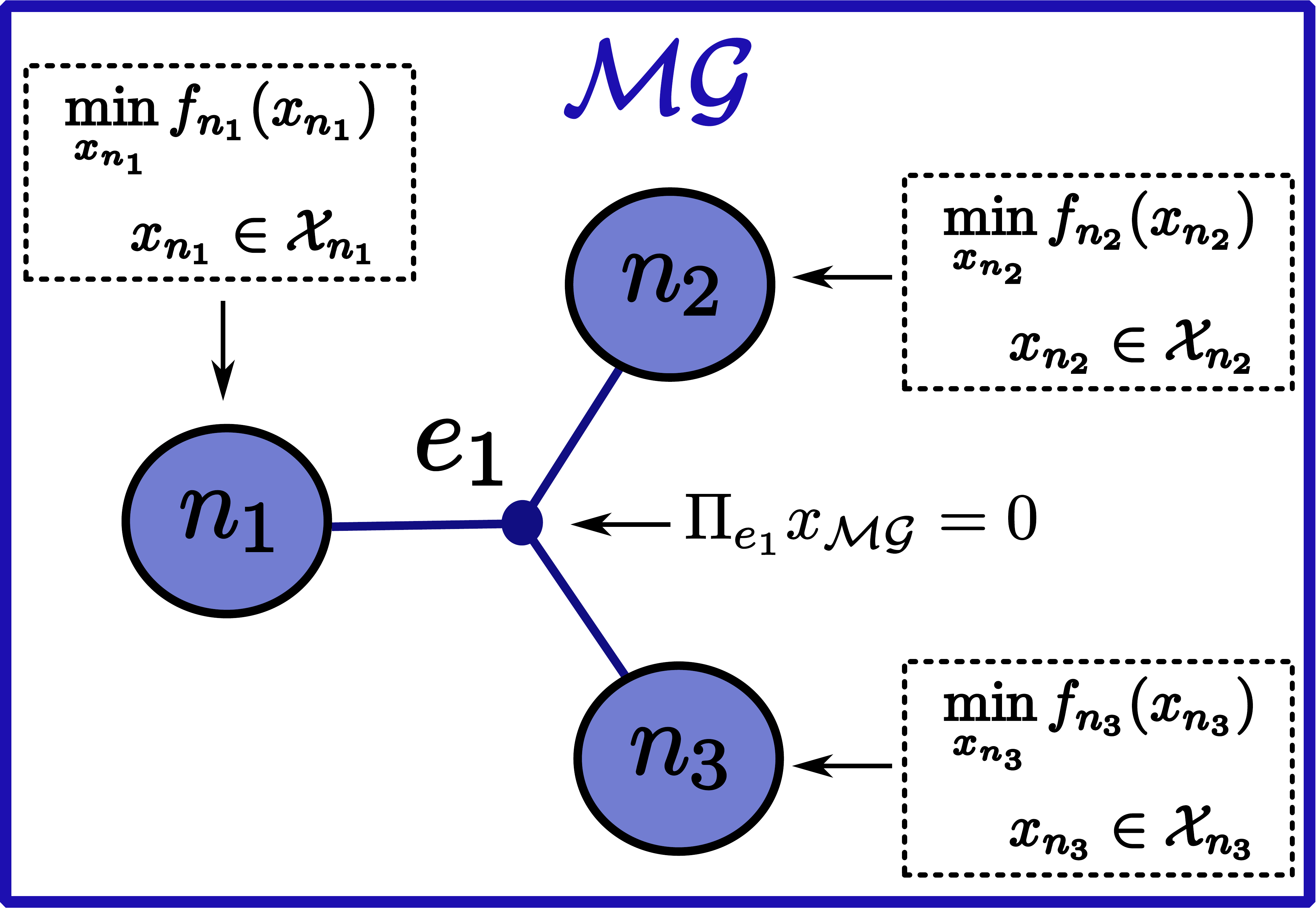}
    \caption{Example of model graph with three nodes and one hyperedge connecting all nodes.  }
    \label{fig:simple_modelgraph}
\end{figure}

It is often useful to define row and column partitions of the connectivity matrix $\Pi_{\mathcal{MG}}$ in order to convey structure to optimization solvers.  If we partition by rows (or columns), we recover how each link connects the model graph nodes (or edges). This is illustrated in Figure \ref{fig:simple_modelgraph} and gives rise to the connectivity matrices:
\begin{subequations}
\begin{align}
    \Pi_{e} &:= (\Pi_{e,n_1},\Pi_{e,n_2},...,\Pi_{e,n_{|\mathcal{N}|}}), \; e\in\mathcal{E}(\mathcal{MG})\\
    \Pi_{n} &:= (\Pi_{e_1,n},\Pi_{e_2,n},...,\Pi_{e_{|\mathcal{E}|},n})^T,\; n\in\mathcal{N}(\mathcal{MG}).
\end{align}
\end{subequations}

\begin{figure}[ht]
    \centering
    \includegraphics[scale=0.45]{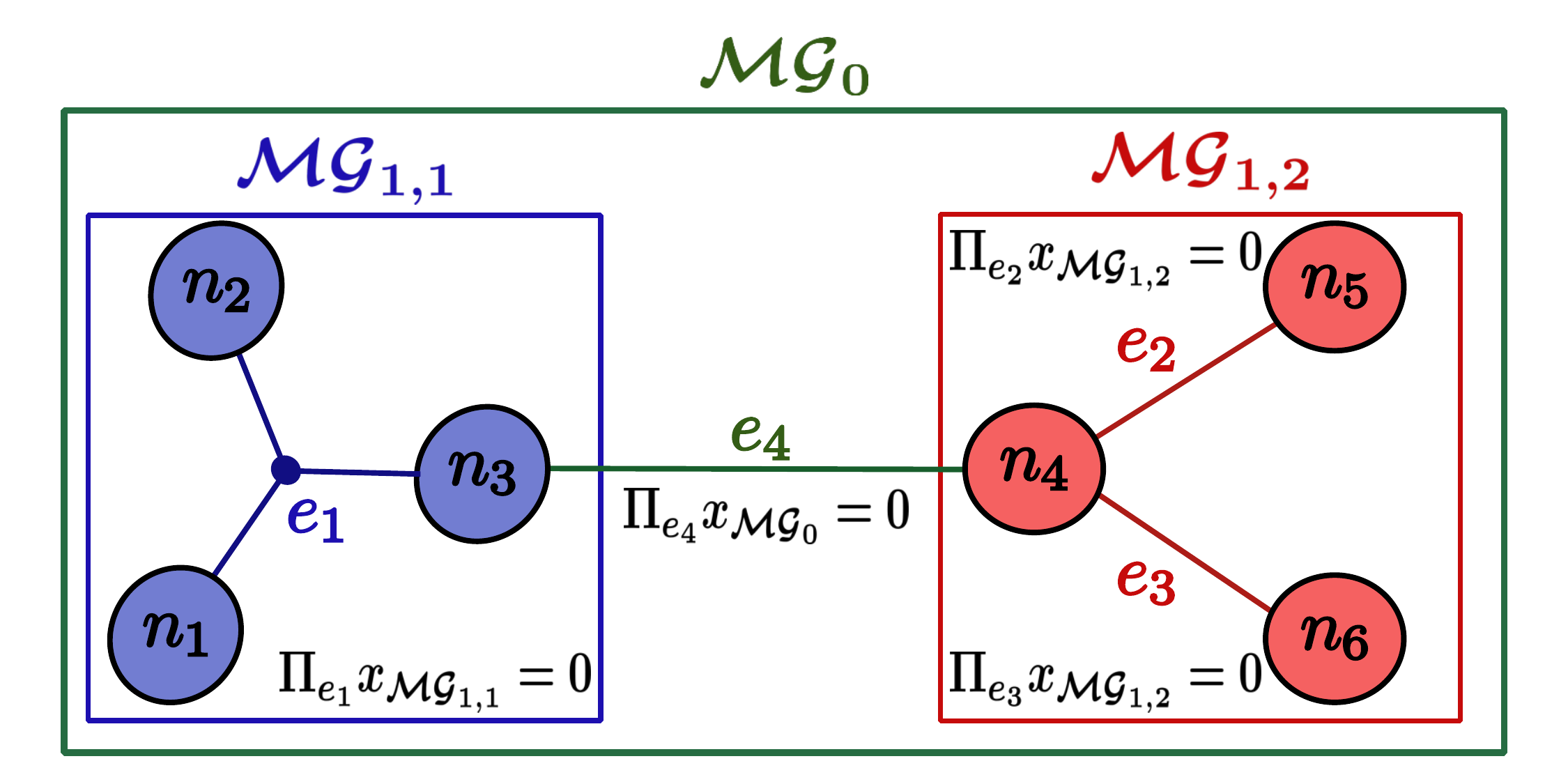}
    \caption{Example of a two-level hierarchical model graph.  $\mathcal{MG}_{1,1}$ contains three nodes ($n_1,n_2,n_3$) and one linking constraint ($e_1$)
    while $\mathcal{MG}_{1,2}$ contains three nodes ($n_4,n_5,n_6$) and two linking constraints (edges $e_2$ and $e_3$).
    $\mathcal{MG}_0$ is the base (top) layer of the graph containing all six nodes and connects the two subgraphs with a single constraint between nodes $n_3$ and $n_4$ (edge $e_4$).}
    \label{fig:subgraph_structures}
\end{figure}

\subsection{Hierarchical Graphs}

The concept of an algebraic graph facilitates the creation of \emph{hierarchical} graphs, wherein a node can represent an algebriac graph
(see Figure \ref{fig:subgraph_structures}). This representation arises in many applications such as integrated planning, scheduling, hierarchical control, coupled power transmission and distribution networks,
and multi-stage stochastic programming. For example, the two-level hierarchical model graph corresponding to Figure \ref{fig:subgraph_structures} can be represented as:

\begin{subequations}\label{eq:model-graph}
\begin{align}
        \min_{x_{\mathcal{MG}_0}} & \quad \sum_{n \in \mathcal{N}(\mathcal{MG}_0)} f_n(x_n) &\textrm{(Graph objective function})\\
        \textrm{s.t.} & \quad x_n \in \mathcal{X}_n, \quad n \in \mathcal{N}(\mathcal{MG}_0) \quad &(\textrm{Local node constraints})\\
        & \Pi_{\mathcal{MG}_0} x_{\mathcal{MG}_0} = 0  \quad &(\textrm{Graph link constraints})\\
        & \Pi_{\mathcal{MG}_{1,1}} x_{\mathcal{MG}_{1,1}} = 0  \quad &(\textrm{Subgraph 1 link constraints})\\
        & \Pi_{\mathcal{MG}_{1,2}} x_{\mathcal{MG}_{1,2}} = 0.  \quad &(\textrm{Subgraph 2 link constraints})
\end{align}
\end{subequations}
Here, we have that $\mathcal{N}(\mathcal{MG}_{1,1}) \subseteq \mathcal{N}(\mathcal{MG}_0)$ and $\mathcal{N}(\mathcal{MG}_{1,2}) \subseteq \mathcal{N}(\mathcal{MG}_0)$ are {\em children} graphs
(partitions or subgraphs) of the {\em parent} graph $\mathcal{MG}_0$. Moreover, we have that
$\mathcal{E}(\mathcal{MG}_{1,1}) \not\subset \mathcal{E}(\mathcal{MG}_0)$ and $\mathcal{E}(\mathcal{MG}_{1,2}) \not\subset \mathcal{E}(\mathcal{MG}_0)$.  In other words, the parent graph $\mathcal{MG}_0$
contains every node in every child subgraph but edges are local to their corresponding subgraphs.  Localizing edges in this form facilitates the expression of hierarchy.

We can generalize formulation \eqref{eq:model-graph} to account for an arbitrary number of subgraphs by defining the set of subgraphs $\mathcal{SG}_1(\mathcal{MG}_0)$. This set contains elements $\mathcal{MG}_{1,1},\mathcal{MG}_{1,2},...$ that are derived by partitioning the parent graph $\mathcal{MG}_0$.  This gives the general two-level hierarchical graph problem:

\begin{subequations}\label{eq:modelgraph-subgraph-compact}
    \begin{align}
            \min_{x_{\mathcal{MG}_0}} & \quad \sum_{n \in \mathcal{N}(\mathcal{MG}_0)} f_n(x_n) & \textrm{(Graph objective function})\\
            \textrm{s.t.} & \quad x_n \in \mathcal{X}_n, \quad n \in \mathcal{N}(\mathcal{MG}_0) \quad &(\textrm{Local node constraints})\\
            & \Pi_{\mathcal{MG}_0} x_{\mathcal{MG}_0} = 0  \quad &(\textrm{Graph links})\\
            & \Pi_{\mathcal{MG}_1} x_{\mathcal{MG}_1} = 0  \quad  \mathcal{MG}_1 \in \mathcal{SG}_1(\mathcal{MG}_0). &(\textrm{Subgraph links})
    \end{align}
\end{subequations}

Every element $\mathcal{MG}_1$ of the subgraph $\mathcal{SG}_1(\mathcal{MG}_0)$ can in turn be partitioned into a set of of
subgraphs $\mathcal{SG}_2(\mathcal{MG}_1)$ (with elements $\mathcal{MG}_{2,1},\mathcal{MG}_{2,2},...$) to create a three-level hierarchical graph of the form:

\begin{subequations}\label{eq:modelgraph-subsubgraph-compact}
    \begin{align}
            \min_{x_{\mathcal{MG}_0}} & \quad \sum_{n \in \mathcal{N}(\mathcal{MG}_0)} f_n(x_n) & \textrm{(Graph objective function})\\
            \textrm{s.t.} & \quad x_n \in \mathcal{X}_n, \quad n \in \mathcal{N}(\mathcal{MG}_0) \quad &(\textrm{Local node constraints})\\
            & \Pi_{\mathcal{MG}_0} x_{\mathcal{MG}_0} = 0  \quad &(\textrm{Graph links})\\
            & \Pi_{\mathcal{MG}_1} x_{\mathcal{MG}_1} = 0  \quad  \mathcal{MG}_1 \in \mathcal{SG}_1(\mathcal{MG}_0) &(\textrm{Subgraph links})\\
                        & \Pi_{\mathcal{MG}_2} x_{\mathcal{MG}_2} = 0  \quad  \mathcal{MG}_2 \in \mathcal{SG}_2(\mathcal{MG}_1), \mathcal{MG}_1 \in \mathcal{SG}_1(\mathcal{MG}_0). &(\textrm{Subsubgraph links})
    \end{align}
\end{subequations}
The partitioning procedure can be repeated to create an arbitrary number of levels. As can be seen, the graph representation can be conveniently used to arrange complex models.

\subsection{Decomposition Strategies}

Graph models directly provide structural information that facilitates the implementation of decomposition algorithms. For example, formulation \eqref{eq:model-graph-compact} has a partially separable structure, because eliminating the linking constraints \eqref{eq:link-constraint}
results in a fully separable problem
(i.e., each component model $n \in \mathcal{N}(\mathcal{MG})$ can be solved independently).  Continuous problems with partially separable structure induce specialized linear algebra kernels that can be exploited inside solvers, while structures in mixed-integer problems can be exploited using techniques such as Lagrangian decomposition.  For instance, a continuous variant of problem \eqref{eq:model-graph-compact}  with feasible sets of the form $\mathcal{X}_n:=\{x\,|\,c_n(x)=0\}$ gives rise to the following Karush-Kuhn-Tucker (KKT) system:
\begin{subequations}\label{eq:graph-KKT}
    \begin{align}
        & \sum_{n \in \mathcal{N}(\mathcal{MG})} \Big( \nabla_{x_n} f_n(x_n) + \nabla_{x_n} c_n(x_n) \lambda_n \Big ) + \Pi_{\mathcal{MG}} \lambda_{\mathcal{MG}} = 0\\
        & c_n(x_n) = 0, \quad  n \in \mathcal{N}(\mathcal{MG})\\
        & \Pi_{\mathcal{MG}} x_\mathcal{MG} = 0.
    \end{align}
\end{subequations}
Upon linearization, this system of algebraic equations gives rise to the block-bordered  system:
\begin{align}\label{eq:schur}
\left[\begin{array}{cccc|c}
K_{n_1}& &&&\Pi_{n_1}^T\\
&K_{n_2}& &&\Pi_{n_2}^T\\
&&\ddots&&\vdots\\
&&&K_{n_{|\mathcal{N}|}}&\Pi_{n_{|\mathcal{N}|}}^T\\\hline
\Pi_{n_1}&\Pi_{n_2}&\hdots&\Pi_{n_{|\mathcal{N}|}}& \\
\end{array}\right]
\left[\begin{array}{c}\Delta w_{n_1}\\ \Delta w_{n_2}\\ \vdots \\ \Delta w_{n_{|\mathcal{N}|}}\\ \hline \Delta \lambda_{\mathcal{MG}} \end{array}\right]=
-\left[\begin{array}{c}\nabla_{w_{n_1}}\mathcal{L}\\ \nabla_{w_{n_2}}\mathcal{L}\\  \vdots \\ \nabla_{w_{n_{|\mathcal{N}|}}}\mathcal{L}\\ \hline \ \nabla_{\lambda_{\mathcal{MG}}}\mathcal{L}\ \end{array}\right],
\end{align}
where $\mathcal{L}$ is the Lagrange function, $\Delta w_n := (\Delta x_n, \Delta \lambda_n)$ is the primal-dual step and
\begin{align}\label{eq:node_blick}
K_n := \left[\begin{array}{cc}
W_n& J_n^T\\
J_n &
\end{array}\right],
\end{align}
is a block matrix corresponding to node $n$, where $W_n$ is the Hessian of the Lagrange function, $J_n := \nabla_{x_n}c_n(x_n)^T$ is the constraint Jacobian, and $g_n := \nabla_{x_n} f_n(x_n)$ is the objective gradient. It is well-known that the block-bordered structure can be exploited by using Schur decomposition or block preconditioning strategies \cite{pipsnlp,kang2015nonlinear}.

The algebraic graph abstraction facilitates the implementation of distributed algorithms such as Benders, ADMM, or Lagrangian decomposition\cite{Grossmann2012,conejo2006,rodriguez2018benchmarking}.
For instance, in a Lagrangian decomposition scheme applied to \eqref{eq:model-graph-compact}, one solves the node subproblems:
    \begin{align}\label{eq:subprobs}
    {\min}_{x_n\in\mathcal{X}_n} & \quad f_n(x_n)+\lambda_{\mathcal{MG}}^T\Pi_{\mathcal{MG}}x_{\mathcal{MG}},
    \end{align}
in parallel for all $n\in\mathcal{MG}$ and for fixed values of the dual variables $\lambda_{\mathcal{MG}}$. A coordination step is then performed by updating the multipliers of the linking constraints
as $\lambda_{\mathcal{MG}}\leftarrow \lambda_{\mathcal{MG}}+\Pi_{\mathcal{MG}}x_{\mathcal{MG}}$.
In a Benders scheme, it is assumed that the graph has the structure of a tree (variables in a parent node are connected to those in the children nodes but no connectivity is present among children).
In this scheme, multiplier information from the node subproblems is used to construct a master problem that updates the coupling variables.

It is also possible to apply different graph analysis techniques directly to the model graph topology to perform functions such as graph partitioning \cite{Karypis1998},
community detection \cite{FORTUNATO201075}, and
identification of spanning trees.  These strategies can be used, for instance,
to create subgraphs that share a minimum set of linking constraints. Such strategies have shown to improve the performance of decomposition algorithms \cite{Tang2017}.

\subsection{Data Management}

Besides decomposition, there are several other key advantages of using a model graph abstraction. In particular, component models are isolated from the graph topology.
A benefit of this is that it is possible to apply automatic differentiation or convexification techniques (or other processing techniques) to each component model $n \in \mathcal{N}(\mathcal{MG})$ separately,
which often results in computational savings.  It is also possible to exchange components models in a graph without altering the core topology.
Moreover, model syntax remains local to the node and it is thus possible to {\em reuse} a component model template in multiple nodes without having to alter syntax (this is not easy to do in algebraic modeling languages such as {\tt AMPL} or {\tt GAMS}). This feature enables modularity and re-usability and enables the implementation of models as parametric functions of \emph{data}. To highlight this, we consider the following model graph representation:
\begin{subequations}\label{eq:model-graph-data}
    \begin{align}
        \min_{x_{\mathcal{MG}}} & \quad \sum_{n \in \mathcal{N}(\mathcal{MG})} f_n(x_n,\eta_n) \\
        \textrm{s.t.} & \quad x_n \in \mathcal{X}_n(\eta_n) \\
        & \quad \Pi_{\mathcal{MG}} x_{\mathcal{MG}} = 0,
    \end{align}
\end{subequations}
where $\eta_n$ is an input data (attribute) vector associated with the model of node $n$. This modularization approach facilitates the implementation of warm-starting procedures, 
which is key in control or estimation applications  \cite{jalving_optimization-based_2018}.
We also note that the output data (solution attributes) structure inherit the structure of the graph model, thus facilitating analysis and post-processing.

\subsection{Model Graph ${\tt Plasmo.jl}$ Implementation}

We describe the basic elements of a software implementation of the algebraic graph abstraction in ${\tt Plasmo.jl}$ (\url{https://github.com/zavalab/Plasmo.jl}).
${\tt Plasmo.jl}$ is written in the Julia programming language and leverages basic algebraic modeling capabilities of JuMP \cite{DunningHuchetteLubin2017} to express and process individual
node models and the graph analysis capabilities of
${\tt LightGraphs.jl}$ \cite{Bromberger17} to construct and manage graph structures.

A ${\tt Plasmo.jl}$ implementation corresponding to the system shown in Figure \ref{fig:simple_modelgraph} is presented in Figure \ref{fig:model_graphcode_snippet}.  The code snippet shows
how to construct a graph \emph{object} and how to populate it with nodes in lines \ref{line:mg_start} - \ref{line:mg_end}.
Notice that nodes do not require component models upon creation and remain separate from such models to
facilitate modification of the component models (facilitating model swapping and warm-starts).
Next, component models are created and added to each node in lines \ref{line:m1}, \ref{line:m2}, and \ref{line:m3}.
In this example, all nodes contain a local variable named $x$ which are linked together using a model graph link constraint on
line \ref{line:linkconstraint}. Note, again, that syntax of a component model remains private. Finally, the model graph object is mapped into an object that can be solved by an optimization solver (line \ref{line:solve}). The solution of each component model is queried individually on lines \ref{line:sol1} and \ref{line:sol2} (the solution has the same structure as the graph).

\begin{figure}[h]
%\centering
\begin{scriptsize}
\lstset{language=Julia,breaklines = true,xleftmargin=\parindent}
\begin{lstlisting}[escapeinside={(*}{*)},escapechar = |]
     #load Ipopt
     using Ipopt
     #load Plasmo.jl
     using Plasmo
     #create model graph
     mg = ModelGraph(solver = IpoptSolver()) |\label{line:mg_start}|
     #add nodes to graph
     n1 = addnode!(mg)
     n2 = addnode!(mg)
     n3 = addnode!(mg)                        |\label{line:mg_end}|
     #add models to nodes
     setmodel(n1,node_model_1()) |\label{line:m1}|
     setmodel(n2,node_model_2()) |\label{line:m2}|
     setmodel(n3,node_model_3()) |\label{line:m3}|
     #create the model graph link constraint
     @linkconstraint(mg,n1[:x] + n2[:x] + n3[:x] == 0) |\label{line:linkconstraint}|
     # solve model graph
     solve(mg) |\label{line:solve}|
     #query nodes for their solutions
     getvalue(n1[:x])  |\label{line:sol1}|
     getvalue(n2[:x])  |\label{line:sol2}|
\end{lstlisting}
\end{scriptsize} %
\caption{Code snippet for constructing and solving a simple model graph}
\label{fig:model_graphcode_snippet}
\end{figure}

\FloatBarrier

\begin{figure}[!htb]
    \centering
    \includegraphics[scale=0.4]{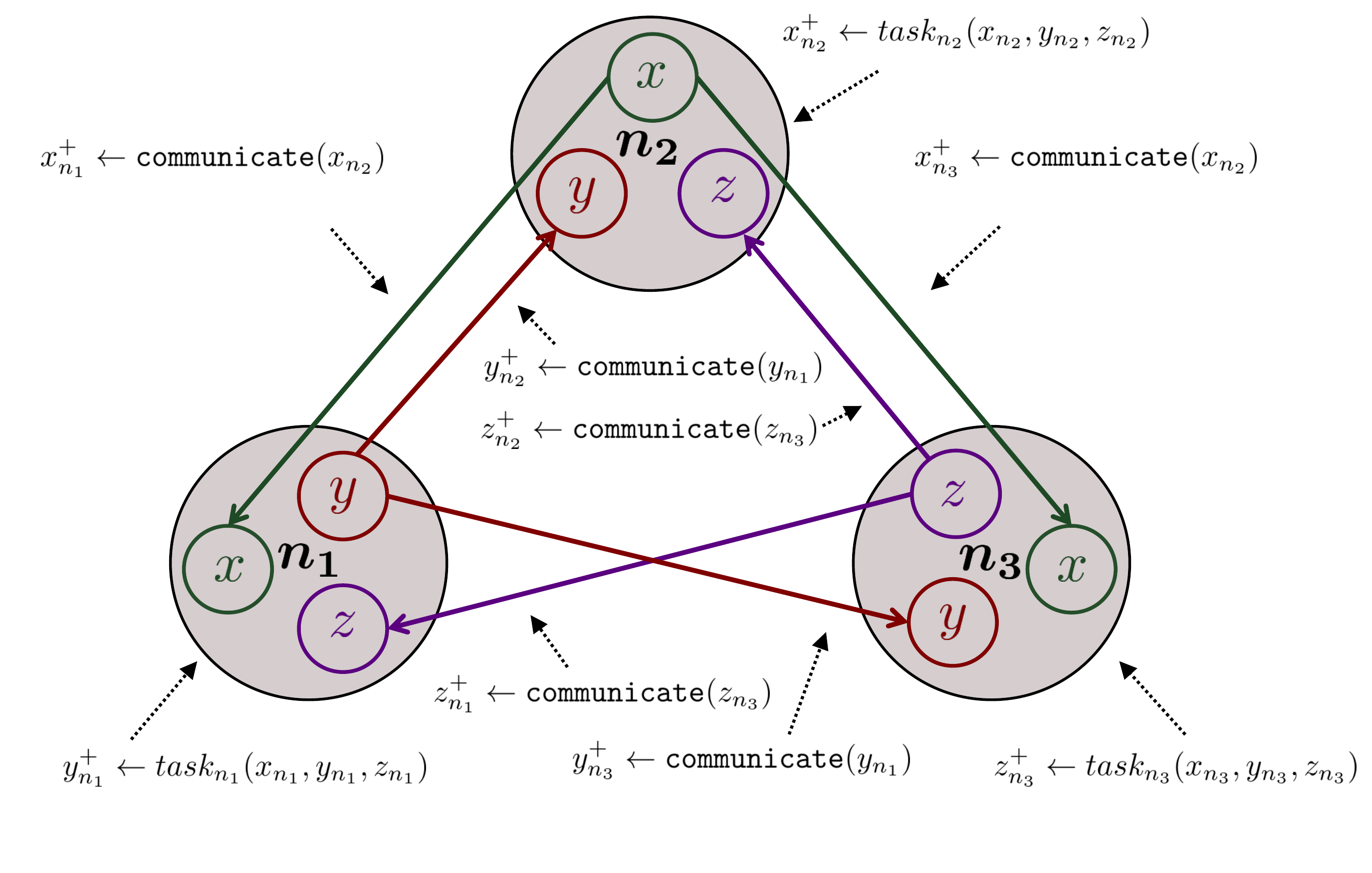}
    \caption{Depiction of a computing graph with three nodes and six edges.  Node $n_1$ computes $task_{n_1}$ using the data attributes ($x,y,$ and $z$) and updates the value of attribute $y$.
    Similarly, node $n_2$ computes $task_{n_2}$ and updates attribute $x$, and node $n_3$ computes $task_{n_3}$ and updates attribute $z$.  Attribute values are
    communicated between nodes using edges.}
    \label{fig:workflows}
\end{figure}

\section{Computing Graphs}\label{sec:workflows}

Simulating cyber systems requires capturing dynamic and logical aspects that arise in real-time decision-making such as
as time delays, computing/processing latency and failures, and asynchronicity.  For example, we might be interested in predicting how a given control system
will perform on an architecture that contains sensors and computing devices with limited processing capabilities (and thus ling delays) or how a
distributed optimization algorithm will perform on a distributed-memory computing cluster compared to a single central processing unit (CPU).
Communication aspects are particularly challenging to capture in cyber systems, as they often involve complex topologies and frequencies.
This section presents the basic elements of a computing graph abstraction seeks to facilitate modeling of cyber systems.

\subsection{Representation}

In the proposed abstraction, a {computing graph} is a {\em directed} multigraph that we denote as $\mathcal{CG}(\mathcal{N},\mathcal{E})$ and that contains a set of nodes $\mathcal{N}(\mathcal{CG})$ (which
perform computing tasks) and edges $\mathcal{E}(\mathcal{CG})$ (which communicate attributes between nodes). 

A {node} $n \in \mathcal{N}(\mathcal{CG})$ contains a set of {\em attributes} $\mathcal{A}_n$ and computing {\em tasks} $\mathcal{T}_n$.
The attributes $\mathcal{A}_n$ represent node data and tasks $\mathcal{T}_n$ are computational procedures that operate on and/or change attributes.
In other words, {\em a computing task maps attributes} (a task takes attribute data and processes it to create other attribute data).
This interpretation resembles that of a manufacturing process, which takes raw material to generate products. Each task $t\in\mathcal{T}_n$ in the computing graph requires a given execution time $\Delta \theta_{t}$. An {edge} $e \in \mathcal{E}(\mathcal{CG})$ contains a set of
attributes $\mathcal{A}_e$ associated with its support nodes $\mathcal{N}(e)$. Communicating attributes between nodes involves a communication delay $\Delta \theta_{e}$.

The collection of computing and communication tasks comprises an {\em algorithm} (also known as a computing workflow in the computer science community).
Consequently, a computing graph abstraction seeks to facilitate the {\em creation and simulation of algorithms}.

The computing graph contains
a {\em global clock} $\theta$ and each node has an internal {\em local clock} $\theta_n$. The clocks are used to manage and schedule computing tasks and communication. For any task $t$ executed at time $\theta$,
its attributes become updated at clock time $\theta + \Delta \theta_{t}$.
Likewise, for any edge $e$ that communicates its attribute at time $\theta$,
the destination attribute value is updated with the source attribute value at time $\theta  + \Delta \theta_e$.
Under the proposed abstraction, computing and communication tasks can be {\em synchronous} (a task is not executed until all attributes are received) or {\em asynchronous} (a task is executed with current values).
This enables capturing a wide range of behaviors seen in applications.

The computing and communication times $\Delta \theta_t$ and $\Delta \theta_e$ can represent {\em true}
times (times required by the computing devices executing the tasks) or {\em virtual} times (times required by hypothetical devices executing tasks).
In other words, the proposed abstraction allows us to simulate the behavior of algorithms on virtual (hypothetical) computing architectures.
This is beneficial when one lacks access to an actual sophisticated architecture (such as a large-scale parallel computer or an industrial control system) but one seeks to predict how effective an algorithm will be when executed under such architecture. Moreover, this enables us to analyze the behavior of algorithms under extreme events that might involve communication or computation failures and with this test their resilience.

Both nodes and edges use the concept of {\em state managers} to manage task behavior (e.g., determining when a task has been completed) and to manage
communication (e.g., determining when data is sent or received). This representation resembles that used in process scheduling and has interesting connections with automata theory
and discrete event simulation \cite{agarwal2009linear}. These connections can be exploited to derive a coherent state-space representation.
In Section \ref{app:compute-graph} we provide details on how this can be used to facilitate implementation of the computing graph abstraction.

Figure \ref{fig:workflows} depicts an example computing graph containing three nodes and six edges.  Each node contains a single task which takes local attribute values
$x$, $y$, and $z$ as input and updates one of their values.  For example, $task_{n_1}$ processes its attributes and updates the value of attribute $y$.
The nodes communicate attribute values with each other using the six edges. For instance, attribute $y$ is communicated to both nodes $n_2$ and $n_3$ which updates
the value of $y$ on these respective nodes. The superscript $+$ indicates that the attribute may be updated after a given time
(to capture computing and communication delays).

It is important to highlight the differences and synergies between a computing and an algebraic graph. In a computing graph a node contains a dynamic component (a computing task) while in a algebraic graph a
node contains a static component (an algebraic model). Moreover, in a computing graph an edge connects attributes (dynamically) while in a model graph an edge connects algebraic variables (statically).
Under a computing graph, the {\em solution of an algebraic graph} is considered a computing task. Consequently, a computing node might use an algebraic graph to perform a given task or a computing graph might be an algorithm for solving a given algebraic graph. For instance, for the former, we might create a computing graph that executes a control algorithm and use an algebraic graph to simulate the behavior of the physical system under the actions of the control system. For the later, we might create a computing graph that executes  a solution algorithm (e.g.,  Bender decomposition) to solve an algebraic graph. These capabilities enable the simulation of complex cyber-physical systems.  In Section \ref{sec:case_study}, we provide examples on how this can be done.  We also highlight that computing tasks are general and might involve procedures that go beyond the solution of algebraic graphs such as forecasting, data analysis, learning, solution of optimization problems (that are not expressed as graphs), and so on.

\subsection{Computing Graph ${\tt Plasmo.jl}$ Implementation}

Here, we present an implementation example in ${\tt Plasmo.jl}$ to highlight the features of computing graphs. The example corresponds to that of Figure \ref{fig:workflows}
(this structure resembles that of a distributed optimization or control algorithm). The ${\tt Plasmo.jl}$ implementation is shown in Figure \ref{fig:workflow_code_snippet}.
Line \ref{line:workflow} creates the computing graph, line \ref{line:wnode1} adds the first node and line \ref{line:node1attributes}
adds three attributes to the node (named $x$, $y$, and $z$).  We also specify {\tt start = 0} on this line to provide an initial value for the data attributes.
We add a task to node $n_1$ in line \ref{line:task1} which runs the task function
$task_{n1}$ and provide the keyword argument {\tt triggered\_by = Received(x,z)} to indicate that receiving attributes $x$ or $z$ will trigger this task.  We also provide
keyword argument {\tt compute\_time = :walltime}, which indicates that $\Delta \theta_{task_{n_1}}$ is the {\em true} compute time of the function (in the computing device executing it),
and {\tt trigger\_during\_busy = :queue\_task}
which indicates that triggering the task during computation will queue the task until its current computation is finished.

Next, we add the second computing node with the same attributes in lines \ref{line:wnode2} and \ref{line:node2attributes}.  We add the
task $task_{n_2}$ to the node in line \ref{line:task2} and provide similar keyword arguments.  We add the last node with attributes in lines \ref{line:wnode3} and \ref{line:node3attributes}
add its task in line \ref{line:task3}.  For this node, we use the keyword argument {\tt triggered\_by = Updated(z)}, which indicates that updating attribute z will trigger that task (i.e., if $task_{n_3}$ updates
attribute $z$, the task will run continuously).  Here, we fix the computing time to two time units using the keyword argument {\tt compute\_time = 2}. This forces the computing time of the task to be a
fixed value (as opposed to the actual value of the processor). This feature allows us to simulate virtual behavior or a hypothetical computing device.

Node $n_1$ is connected to nodes $n_2$ and $n_3$ in line \ref{line:connect1} which
communicates the value of attribute $y$ on node $n_1$ to attribute $y$ on nodes $n_2$ and $n_3$.  We provide the keyword argument {\tt delay = 2} which indicates that communication requires a fixed time of two time units (to simulate a fixed delay)
and the argument {\tt send\_on = Updated(n1[:y])}, which indicates that communication occurs when the source attribute value {\tt n1[:y]} is updated (attribute $y$ residing in node $n_1$). 
We also connect attribute $x$ in node $n_2$ to that in nodes $n_1$ and $n_3$ in line \ref{line:connect2} and connect attribute $z$ in node $n_3$ to nodes $n_1$ and $n_2$ in line \ref{line:connect3}.

We provide an initial trigger to start each task at global time $\theta = 0$ and execute the
computing graph for twenty time units in lines \ref{line:schedule} and \ref{line:execute}.

\begin{figure}[h!]
\centering
\begin{scriptsize}
\lstset{language=Julia}
\begin{lstlisting}[escapeinside={(*}{*)},escapechar = |]
    using Plasmo
    #create graph
    graph = ComputingGraph()        |\label{line:workflow}|

    #add first node and attributes
    n1 = addnode!(graph)                    |\label{line:wnode1}|
    @attributes(n1,x,y,z,start = 0)           |\label{line:node1attributes}|

    #add task to node
    @nodetask(graph,n1,task_n1,triggered_by = Received(x,z), compute_time = :walltime, trigger_during_busy = :queue_task) |\label{line:task1}|

    #add second node and attributes
    n2 = addnode!(graph)         |\label{line:wnode2}|
    @attributes(n2,x,y,z,start = 0)           |\label{line:node2attributes}|

	#add task to node
    @nodetask(graph,n2,task_n2,triggered_by = Received(y,z),compute_time = :walltime, trigger_during_busy = :queue_task)  |\label{line:task2}|

    #add third node and attributes
    n3 = addnode!(graph)         |\label{line:wnode3}|
    @attributes(n2,x,y,z,start = 1)           |\label{line:node3attributes}|

    #add task to node
    @nodetask(graph,n3,task_n3, triggered_by=Updated(z),compute_time = 2)      |\label{line:task3}|

	#create edge and link attributes
    @connect(graph,n1[:y] => [n2[:y],n3[:y]], delay=2, send_on=Updated(n1[:y]))     |\label{line:connect1}|

    #create edge and link attributes
    @connect(graph,n2[:x] => [n1[:x],n3[:x]], delay=1, send_on=Updated(n2[:x]))     |\label{line:connect2}|

    #create edge and link attributes
    @connect(graph,n3[:z] => [n1[:z],n2[:z]], delay=1, send_on=Updated(n3[:z]))     |\label{line:connect3}|


    #Initialize graph with an execute signal at time zero
    schedule_trigger(graph,n1[:task_n1],n2[:task_n2],n3[:task_n3],time = 0)  |\label{line:schedule}|
    execute!(graph,20)               |\label{line:execute}|
\end{lstlisting}
\end{scriptsize} %
\caption{Code snippet for creating a computing graph in ${\tt Plasmo.jl}$}
\label{fig:workflow_code_snippet}
\end{figure}

The keyword arguments in ${\tt Plasmo.jl}$ can be used to manage how tasks and edges respond to different signals.
This makes it possible to simulate diverse computation and communication behaviors.  For instance, if each task in Figure \ref{fig:workflows} was to solve 
an optimization problem (e.g., an algebraic graph), the computing graph described in Figure \ref{fig:workflow_code_snippet} could simulate the behavior of a distributed optimization algorithm
wherein nodes $n_1$ and $n_2$ run their optimization problems in response to receiving attribute values and $n_3$ would compute its optimization problem continuously (or asynchronously).

The computing behavior of Figure \ref{fig:workflow_code_snippet} could also be changed depending on the nature of the tasks and arguments provided.  For instance, it would be possible to
simulate a control system wherein one node could continuously compute a simulation task (e.g., a plant simulation) and its outgoing edges could communicate on a sampling interval
(as opposed to communicating based on source updates). 
The other nodes could receive measurement attributes which would trigger their tasks to compute control actions and update their inject attributes which would be communicated back to the plant node.

\subsection{State-Space Representation of Computing Graph}\label{app:compute-graph}

In the proposed abstraction, the computation and communication logic provided by the computing graph is driven by an underlying state manager abstraction wherein 
transitions in tasks states are triggered by input signals. This representation conveniently captures computing and communication latency and is flexible and extensible.  In this section, we provide details on how state managers can facilitate the implementation of computing graphs. We highlight, however, that an actual user does not need to know these details to use {\tt Plasmo.jl}.

\begin{figure}[!htb]
    \centering
    \includegraphics[scale=0.5]{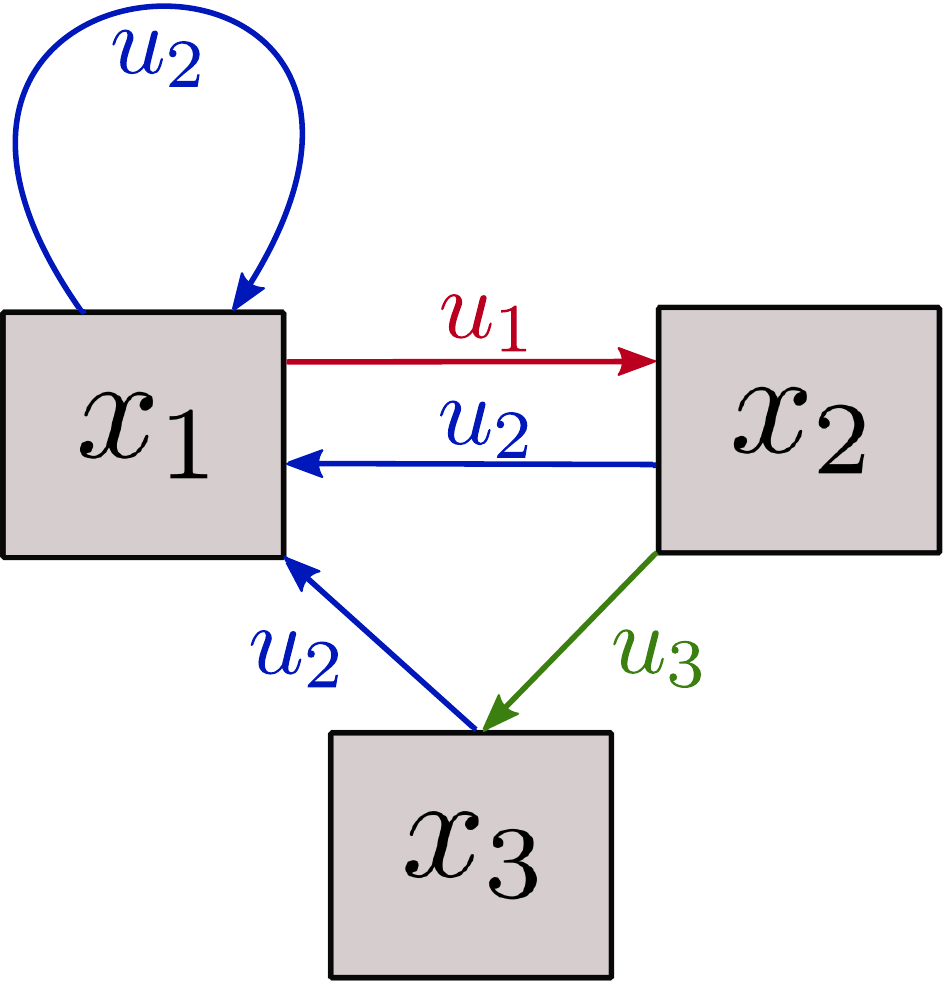}
    \caption{A simple state machine with three states ($x_1,x_2,x_3$), three action signals ($u_1,u_2,u_3$), and five possible state transitions.}
    \label{fig:simple_state_machine}
\end{figure}

Node and edge managers use \emph{states} and \emph{signals} to
manage task computation and attribute communication.  The use of managers is motivated by state machine abstractions from automata theory \cite{Brand1983,Bollig2017},
which are often used in control-logic applications \cite{Endsley2004}.  State machines are also used to manage logical behaviors in Simulink \cite{Branicky1997} and
agent-based simulation frameworks \cite{Ozik1994}.  State machines can be used to represent actions that trigger transitions in task states.  Figure \ref{fig:simple_state_machine} presents a simple state machine
with three states, three signals, and five possible transitions between states. 

In a computing graph, a node manager $M_n$ oversees the state of node tasks. Such tasks are specified by the user in the form of {\em functions}. Each node has an associated tuple ($\mathcal{X}_n,\mathcal{U}_n,\mathcal{Y}_n,f_n,g_n$)
where $\mathcal{X}_n$ is the set of states, $\mathcal{U}_n$ is the set of action signals,
$\mathcal{Y}_n \subseteq \mathcal{A}_n$ is the set of broadcast (output) attributes,
$f_n : \mathcal{X}_n \times \mathcal{U}_n \rightarrow \mathcal{X}_n$ is the state transition mapping,
and $g_n:\mathcal{Y}_n \to \mathcal{Y}_n$ is the attribute update mapping. The node dynamic evolution is represented as a system of the form:
\begin{subequations}\label{eq:state_manager}
    \begin{align}
        x^+_n &= f_n(x_n,u_n)\\
        \eta_n^+&=g_n(\eta_n).
    \end{align}
    \end{subequations}
Here, the next state $x^+_n=f_n(x_n,u_n)$  is the result of the mapping from
the current state $x_n$ and a given action signal $u_n$.
Every state transition $(u_n,x_n)\to x_n^+$ also triggers a transition in the attribute values $\eta_n\to \eta_n^+$ (i.e., tasks update attributes).
Action signals can also be sent to the state managers of other nodes in the form attributes.

The proposed abstraction can incorporate an arbitrary number states and actions but here we provide an example of a basic set of states and actions that can be considered in an actual implementation:
\begin{subequations}\label{eq:node_sets}
    \begin{align}
        \mathcal{X}_n &:= ({\tt idle},{\tt executing\_task}, {\tt finalized\_task})\\
        \mathcal{U}_n &:= ({\tt execute\_task}, {\tt finalize\_task}, {\tt attribute\_updated}, {\tt attribute\_received}).
    \end{align}
\end{subequations}
At a given point in time, a node manager can be in one of the states in $\mathcal{X}_n$.
The set of signals recognized by the manager are $\mathcal{U}_n$ and these will trigger (depending on the current state) a transition between states.
Such signals include, for instance, ${\tt execute\_task}$ or ${\tt attribute\_received}$.  The set of broadcast targets (that receive created signals) are the node
itself and all of its outgoing edges $\mathcal{E}(n)$.  Hence, a node can send signals (in the form of attributes) to itself or its outgoing edges.

Using the sets defined in \eqref{eq:node_sets}, we can define a transition mapping $f_n(\cdot)$
as:
\begin{subequations}\label{eq:node_transition_mapping}
    \begin{align}
        &{\tt executing\_task} \leftarrow ({\tt idle},{\tt execute\_task}) \\
        &{\tt finalized\_task} \leftarrow ({\tt executing\_task},{\tt finalize\_task})\\
        &{\tt idle} \leftarrow ({\tt finalized\_task}, {\tt back\_to\_idle}).
    \end{align}
\end{subequations}
In \eqref{eq:node_transition_mapping} we can see, for instance, that a task transitions to the ${\tt executing\_task}$
state when it receives the corresponding ${\tt execute\_task}$ signal and it transitions to the ${\tt finalized\_task}$ state
when it is executing a task and receives a signal to finalize such task.
The signals to execute or finalize a task are generated by user-defined attributes.
For instance, a user-defined attribute consisting of a flag such as ${\tt convergence}$ or ${\tt max\_iterations}$ can generate a finalize
task signal that in turn triggers a state transition.

An edge manager $M_e$ can be defined in the same way as a node manager with associated states and signals for communication.
A possible implementation of an edge manager includes the following states and actions:
\begin{subequations}\label{eq:edge_sets}
    \begin{align}
        &\mathcal{X}_e := ({\tt idle},{\tt communicating},{\tt all\_received})\\
        &\mathcal{U}_e := ({\tt attribute\_updated}, {\tt communicate}, {\tt attribute\_sent}, {\tt attribute\_received}).
    \end{align}
\end{subequations}
An edge can send signals to itself or its supporting nodes $\mathcal{N}(e)$ in the form of its attributes.
Using the sets defined in \eqref{eq:edge_sets}, we can define a transition mapping $f_e(\cdot)$ of the form:

\begin{subequations}\label{eq:edge_transition_mapping}
    \begin{align}
        &{\tt communicating} \leftarrow ({\tt idle},{\tt communicate})\\
        &{\tt all\_received} \leftarrow ({\tt communicating},{\tt all\_received})\\
        &{\tt idle} \leftarrow ({\tt all\_received}, {\tt back\_to\_idle}).
    \end{align}
\end{subequations}
The mappings in \eqref{eq:edge_transition_mapping} closely reflect the node transition mapping in \eqref{eq:node_transition_mapping}.  An edge
transitions to the communicating state when it receives
a ${\tt communicate}$ signal and transitions to the ${\tt all\_received}$ state when it receives
the ${\tt all\_received}$ signal (indicating that all sent attributes were received).

\begin{figure}[!htb]
    \centering
    \includegraphics[width=6 in]{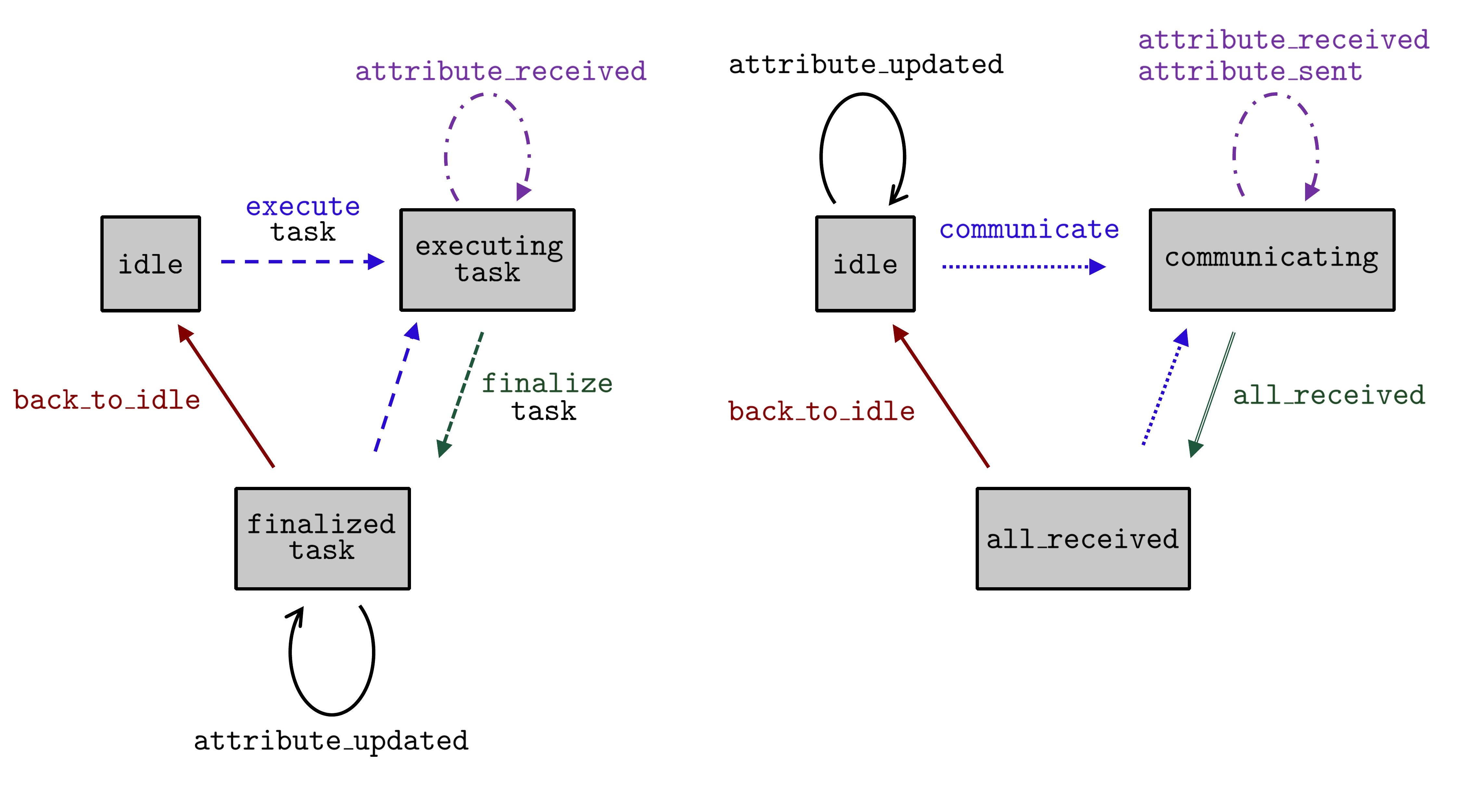}
    \caption{Implementation of a node manager $M_n$ (left) and an edge manager $M_e$ (right).  Action signals
    trigger transitions between states and can include transitions that return to the same state.}
    \label{fig:state_managers}
\end{figure}

Figure \ref{fig:state_managers} depicts the node and edge manager transition mappings \eqref{eq:node_transition_mapping} and \eqref{eq:edge_transition_mapping} with additional transitions.
This figure highlights that action signals can trigger self-transitions wherein nodes or edges loop back to their original state.  For instance, a node can transition
back to the ${\tt executing\_task}$ state when it receives the signal ${\tt attribute\_updated}$ and an edge can transition back to the ${\tt communicating}$ state when
it receives either ${\tt attribute\_sent}$ or ${\tt attribute\_received}$ signals.
Self-transitions allow attribute updates to occur
during task execution or edge communication.

\subsection{Task Scheduling and Timing}

A computing graph implementation needs to capture timings and order of task execution and attribute communication.
These timings can be managed using a discrete-event queue wherein items in the queue are evaluated
in an order based on their scheduled time \cite{Fishman}.  The computing graph specifically uses a
\emph{signal-queue} wherein action signals are evaluated based on their scheduled evaluation time.

\begin{figure}[!htb]
    \centering
    \includegraphics[scale=0.5]{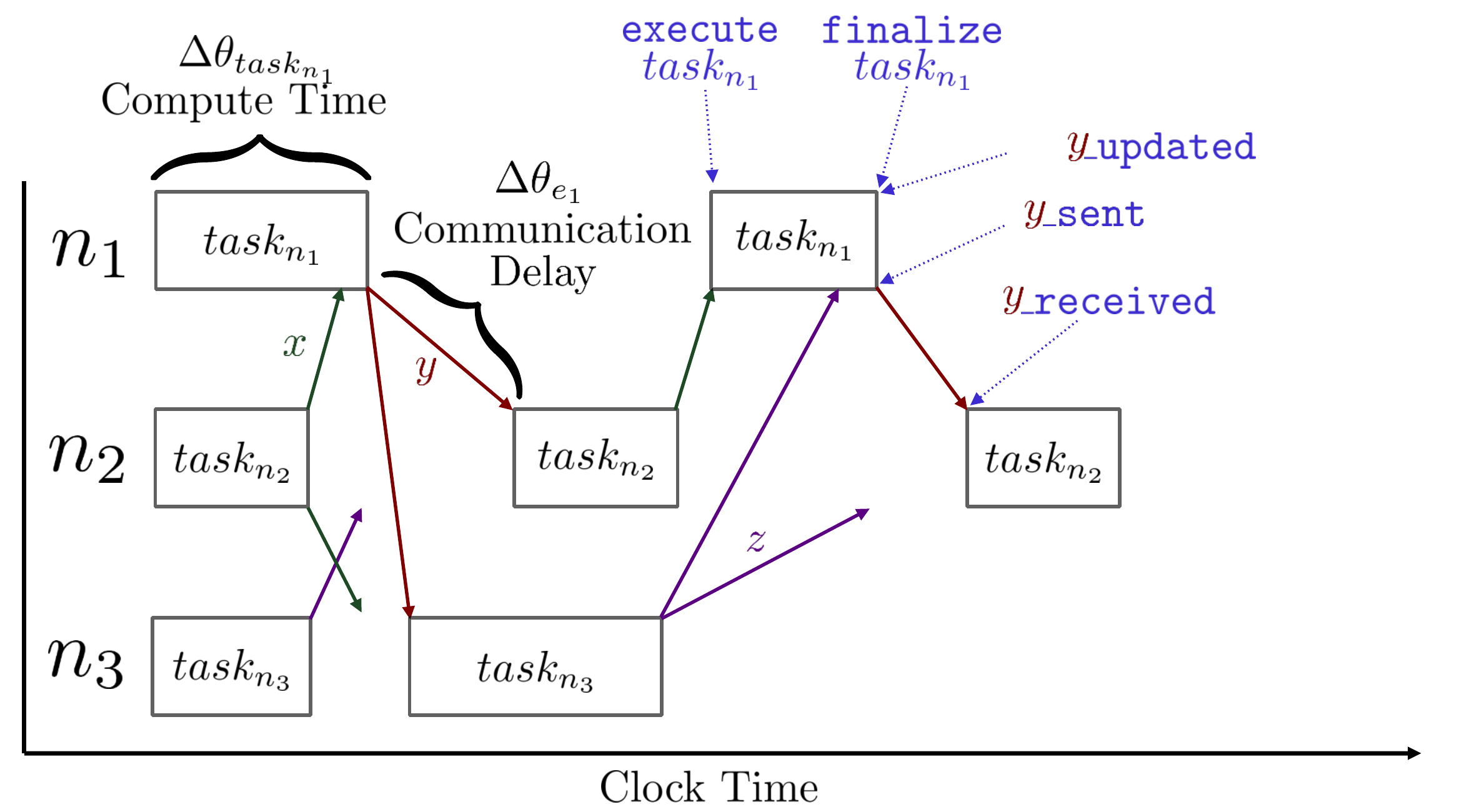}
    \caption{An example execution of the computing graph in Figure \ref{fig:workflows}.  Computing times and communication delays are captured using
    action signals.}
    \label{fig:workflow_queue}
\end{figure}

Figure \ref{fig:workflow_queue} illustrates an example execution of the computing graph from Figure \ref{fig:workflows}.
Node $n_1$ computes $task_{n_1}$ (which requires compute time $\Delta \theta_{task_{n_1}}$) after which the value of attribute $y$ is sent to the corresponding attribute $y$ on
nodes $n_2$ and $n_3$ (which each requires communication time).
The compute and communication times are expressed using signals.  For instance, when
$task_{n_1}$ completes its execution, it produces a ${\tt finalize} \ task_{n_1}$ signal with a delay $\Delta\theta_{task_{n_1}}$ to capture the computation time.
Equivalently, when edge $e_1$ that connects node $n_1$ to node $n_2$ communicates attribute $y$, it produces the $y\_{\tt received}$ signal with a delay $\Delta \theta_{e_1}$.

\FloatBarrier

\section{Case Studies}\label{sec:case_study}
In this section we illustrate how the proposed graph-based abstractions can be used to simulate
complex cyber-physical systems.  We first present a case study which models a pipeline gas network and show how the model graph
facilitates partitioning of the network to enable a scalable solution. We then illustrate how the computing graph can be used to simulate
the behavior of optimization and control algorithms under different computing architectures.

\subsection{Optimal Control of Gas Pipeline Network}
\begin{figure}[!htb]
    \centering
    \includegraphics[scale=0.6]{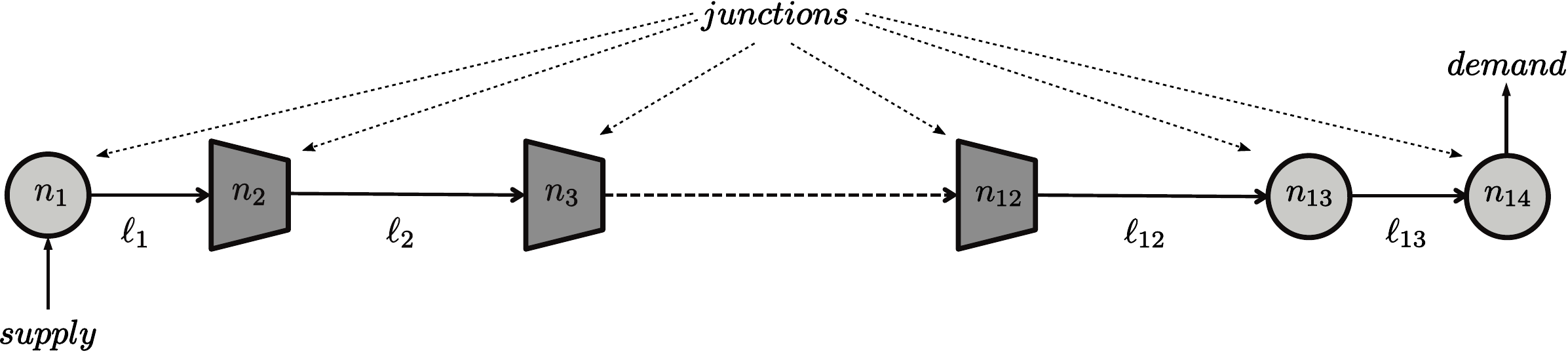}
    \vspace{0.2in}

        \includegraphics[scale=0.35]{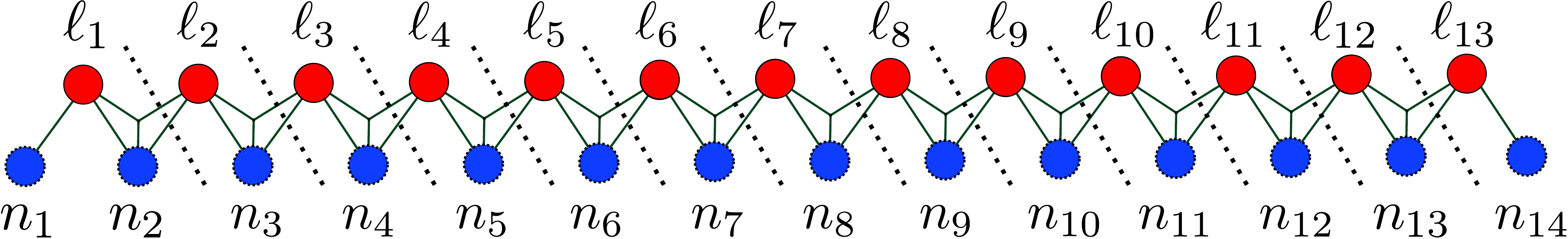}
    \caption{Depiction of pipeline network system (top) and model graph structure partitioned with a k-way scheme (bottom).  Blue nodes represent junction models, red nodes represent pipeline models, and edges indicate
    linking constraints between the models.}
    \label{fig:13_pipeline_model_structure}
\end{figure}

We consider a system of connected pipelines in series \cite{Zavala2014} shown in Figure \ref{fig:13_pipeline_model_structure}.
This linear network includes a gas supply at one end, a time-varying demand at the other end, and twelve compressor stations (giving a total of 14 gas junctions).
The gas junctions connect thirteen pipelines. This forms a model graph with a linear topology. We formulate an optimal control problem to track a sudden demand withdrawal.

Under a model graph abstraction $\mathcal{MG}$, all components (gas junctions and pipelines) are treated as nodes $\mathcal{N}(\mathcal{MG})$ (there are a total of 27 nodes in the model graph).
Note that a pipeline is treated as a node as well and has a corresponding component model (discretized partial differential equations that capture conservation of mass and momentum).
Edges indicate algebraic connectivity between the nodes. The structure of the model graph is depicted in Figure \ref{fig:13_pipeline_model_structure}.
Here, we also show the optimal partitioning of the graph using a $k$-way algorithm.

We solve the optimal control problem by: (i) aggregating all nodes and edges to produce a large optimization problem which ignores the graph structure, and (ii) by exploiting $k$-way partitioning to enable a
parallel Schur decomposition scheme.  The code snippet shown in Figure \ref{fig:code_snippet_13_pipeline} depicts how straightforward it is to implement both approaches using the model graph and highlights
several features that are enabled by using a graph abstraction.  Lines \ref{line:13pipe1} through \ref{line:13pipe2} create the model graph, add nodes, set component models, and connect nodes using link constraints.
For approach (i), the graph is flattened into a traditional
optimization model and solved with {\tt Ipopt} (line \ref{line:13solve1}).  In approach (ii), a k-way partioning scheme \cite{Karypis1998}
is applied to the model graph to produce thirteen partitions and the partition information is passed to the {\tt PIPS-NLP} solver (line \ref{line:13solve2}).
The resolution of the spatial discretization mesh of the pipeline PDEs is gradually increased to produce
optimization problems that span the range of 100,000 to 2 million variables.
The partitioned problem can be solved with near-linear scaling whereas the solution time of the flattened formulation scales cubically.

\begin{figure}[h!]
\centering
\begin{scriptsize}
\lstset{language=Julia}
\begin{lstlisting}[escapeinside={(*}{*)},escapechar = |]
     using GasModels, Plasmo, Ipopt, MPI

	 # get pipeline network data
     include("gas_network_data.jl")

     mg = ModelGraph()    |\label{line:13pipe1}|
     # set junction models
     for i = 1:14
        addnode!(mg)
        setmodel(node,GasNodeModel(node_data[i]))
     end
     # set pipeline (PDE) models
     for j = 1:13
        addnode!(mg)
        setmodel(node,GasPipelineModel(pipeline_data[j]))
     end
     gas_nodes = collect_nodes(mg)[1:14]
     pipelines = collect_nodes(mg)[15:end]
     #Connect gas nodes based on topology
     for node in nodes
         pipes_in = node_map_in[node]
         pipes_out = node_map_out[node]
         @linkconstraint(mg,[t = time_points], 0 == sum(pipes_in[i][:fout][t] for i = 1:length(pipes_in))
         - sum(pipes_out[i][:fin][t] for i = 1:length(pipes_out)) + node[:total_supplied][t] - node[:total_delivered][t])
     end
     #Connect pipelines based on topology
     for pipe in pipelines
         node_from,node_to = edge_map[pipe]
         @linkconstraint(mg,[t = times],pipe[:pin][t] == node_from[:pressure][t])
         @linkconstraint(mg,[t = times],pipe[:pout][t] == node_to[:pressure][t])
     end                |\label{line:13pipe2}|

     #approach (i)
     setsolver(mg,IpoptSolver(linear_solver = "ma57")) # call ipopt to solve flattened model
     solve(mg)  |\label{line:13solve1}|

     #approach (ii)
     partitions = partition(mg,13,method = :kway)  # obtain k-way partitions of graph
     setsolver(mg,PipsnlpSolver(partitions = partitions)) |\label{line:13solve2}| # call pips-nlp to solve structure model
     solve(mg)
 \end{lstlisting}
 \end{scriptsize} %
\caption{Code snippet for gas pipeline network problem in ${\tt Plasmo.jl}$}
\label{fig:code_snippet_13_pipeline}
\end{figure}

We note that approach (i), although limited in computational efficiency, still benefits from 
having an algebraic graph representation. For instance, the same PDE component model is appended to different pipeline nodes and this enables the modular construction of complex network models. 
Moreover, approach (i) can also benefit from warm-starting and updating component models with different demands, operational constraints, or lower-fidelity PDE representations.

\begin{figure}[!h]
    \centering
    \includegraphics[scale=0.4]{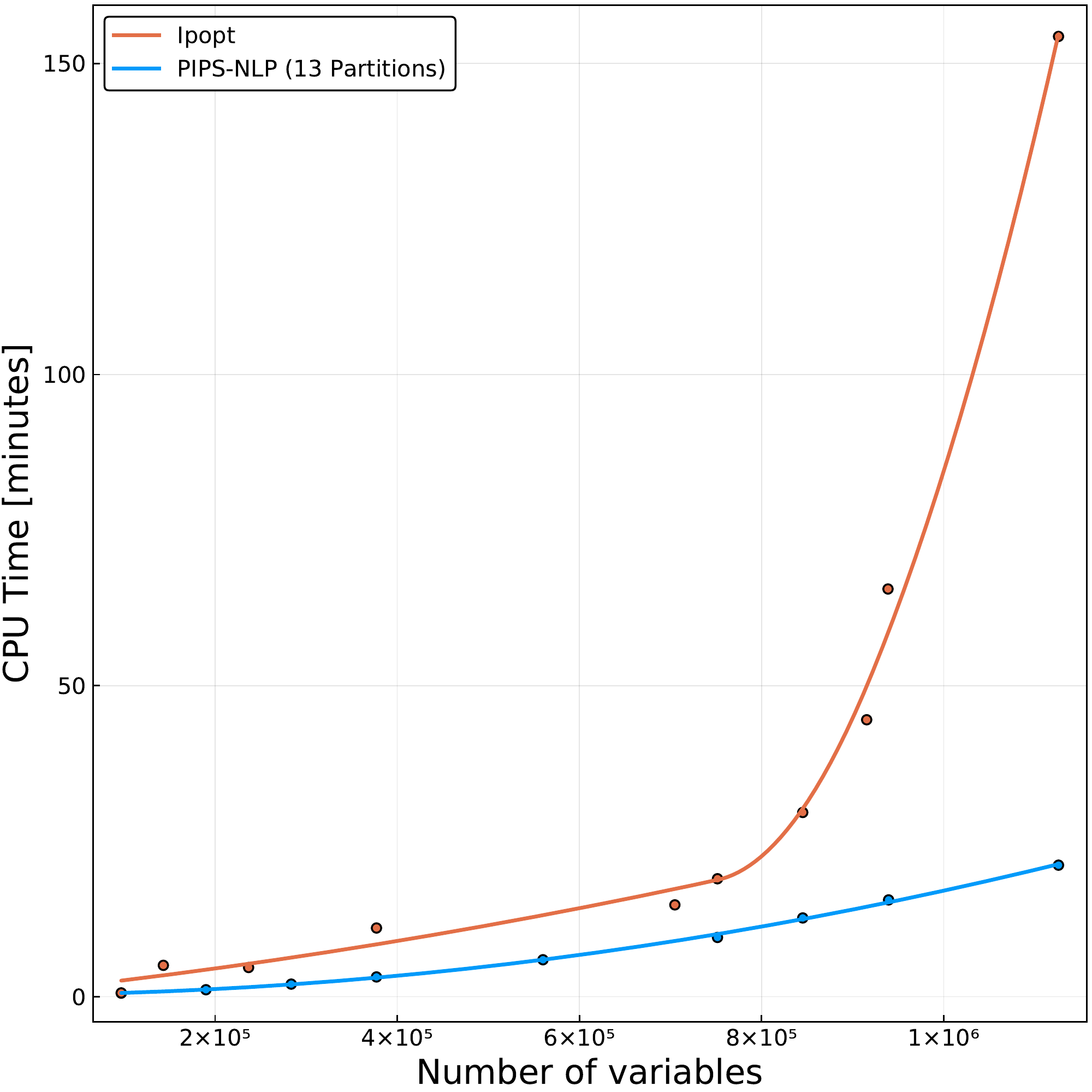}
    \caption{Computational times for the solution of unstructured gas pipeline formulation with {\tt Ipopt} and
    for the solution of the structured formulation using {\tt PIPS-NLP}.}
    \label{fig:13_pipeline_comparison}
\end{figure}

\FloatBarrier

\subsection{Benders Decomposition Algorithm}

We now show how to use computing graphs to simulate the behavior of Benders decomposition under a parallel computing architecture with different numbers of CPUs.
The problem under study is a resource allocation
stochastic program \eqref{eq:benders_problem}.  This stochastic program
can be decomposed with Benders decomposition into a master problem \eqref{eq:benders_master} and a set of subproblems \eqref{eq:benders_sub} defined for a set of sampled scenarios $\Xi$.
This forms a computing graph with a two-level tree topology.

A synchronous variant of Benders decomposition using a master processor and set of $\mathcal{N}:=\{1,...,N\}$ worker CPUs that solve scenario subproblems is described in Algorithm \ref{alg:synchronous_benders}.
The algorithm uses computing tasks implemented in the functions ${\tt solve\_subproblem}$, ${\tt receive\_solution}$, and ${\tt solve\_master}$.
The function ${\tt solve\_subproblem}$ computes a subproblem solution $s$ for a given scenario data sample $\xi\in \Xi$ using the current master solution $\hat{x}$.
Here, the master solution $\hat{x}$, the scenario data $\xi$, and the subproblem solution $s$ are data attributes that are communicated between the master and the worker nodes.
The {\tt solve\_subproblem} function activates the {\tt receive\_solution} task on the master node, which stores subproblem solutions into the set $S$, updates the set of cuts $C$,
and checks how many scenario subproblems have been solved.
If all subproblems are completed, it triggers the execution of the {\tt solve\_master} task using the current set of solution data $S$ and cuts $C$.
Otherwise, it triggers the execution of the {\tt solve\_subproblem} task on a worker processor.  When {\tt solve\_master} is executed, it takes subproblem solution attribute $s$ to update the master
attribute $\hat{x}$ and stops the computing graph if convergence is achieved.  If not converged, it empties the solution set $S$ and activates the {\tt solve\_subproblem} tasks for each worker again,
which will then obtain new subproblem attributes.

\begin{algorithm*}[h!]
\caption{(Synchronous Benders)}
\label{alg:synchronous_benders}
\begin{algorithmic}[1]
    \State \textbf{function} ${\tt solve\_subproblem}$\;($\hat{x}$, $\xi$)
    \State Given the master solution $\hat{x}$ and scenario $\xi$, evaluate sub-problem solution $s$
    \State Activate ${\tt receive\_solution}$($s$) on the master node
    \State \textbf{end function}
    \\\hrulefill
    \State \textbf{function} ${\tt receive\_solution}$($s$)
        \State Given subproblem solution $s$, store solution $S \leftarrow s$, and update cuts $C$
        \If{All subproblems complete}
            \State Activate ${\tt solve\_master}$($S$, $C$) on the master node
        \Else
            \State Get new scenario $\xi$ from scenario set $\Xi$
            \State Activate ${\tt solve\_subproblem}$($\hat{x}$, $\xi$) on worker node
        \EndIf
    \State \textbf{end function}
    \\\hrulefill
    \State \textbf{function} ${\tt solve\_master}$($S$, $C$)
    \State Given solution data $S$ and cuts $C$, update master solution $\hat{x}$
    \State Empty solution data $S \leftarrow \{\}$
    \If{converged}
        \State STOP
    \Else
        \For{$n = 1,...,N$}
        \State Get scenario $\xi_n$ from set $\Xi$
        \State Activate ${\tt solve\_subproblem}$($\hat{x}$, $\xi_n$) on worker node
        \EndFor
    \EndIf
    \State \textbf{end function}
    \\\hrulefill
    \State Initialize variable values and scenario set $\Xi$
    \State Activate ${\tt solve\_master}$ on master node
\end{algorithmic}
\end{algorithm*}

\begin{figure}[!h]
    \centering
    \includegraphics[scale=0.4]{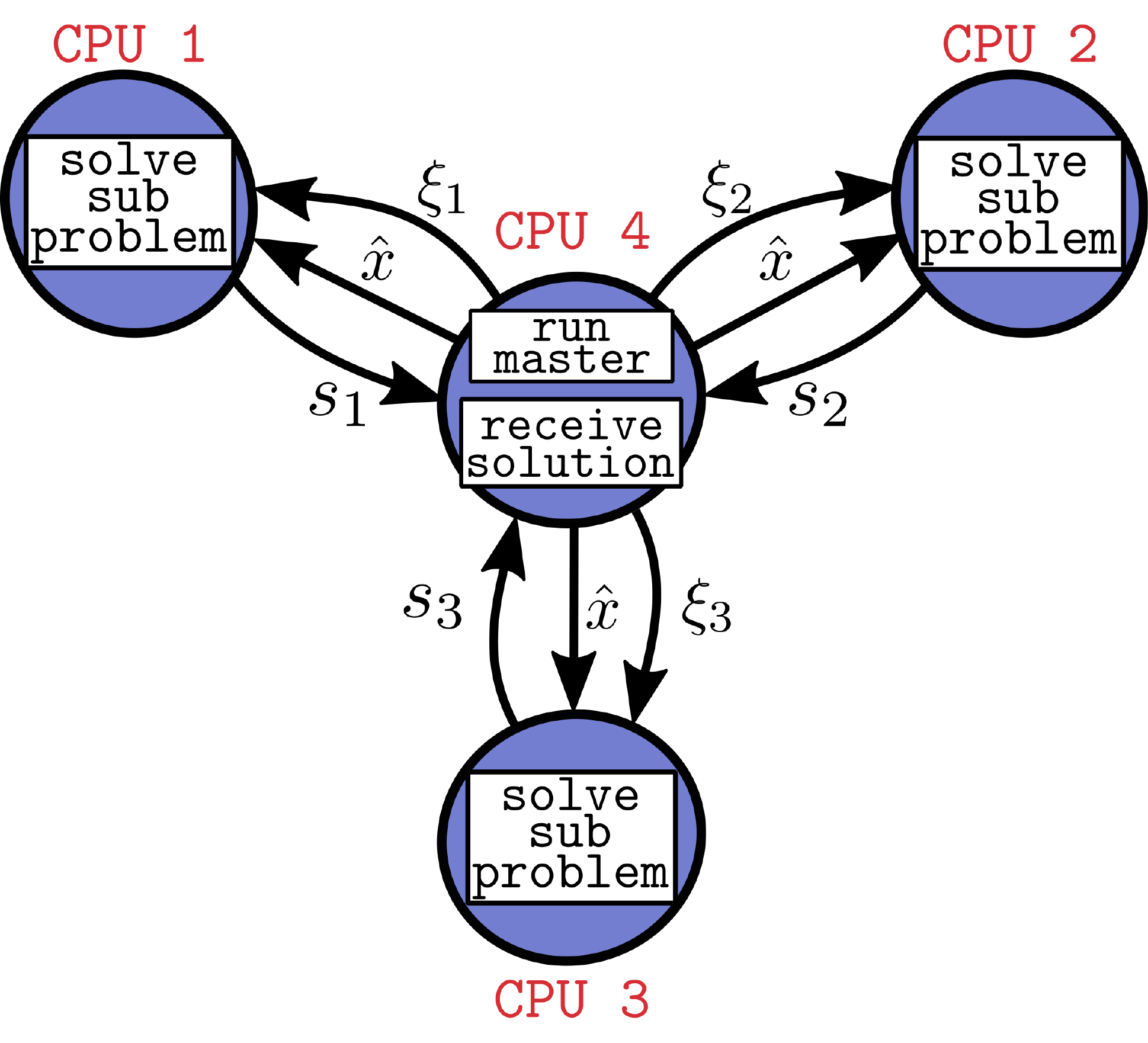}
    \caption{Hypothetical computing architecture executing Benders decomposition.
    CPU 4 executes the solution of the Benders master problem and receives solutions from the subproblems.  CPUs 1, 2, and 3 execute the solution of the scenario subproblems.}
    \label{fig:benders_workflow_topology}
\end{figure}

Algorithm \ref{alg:synchronous_benders} can be modeled and simulated using a computing graph. Under this abstraction, the parent node will solve the master problem and a set of $\mathcal{N}$  children task nodes will
be available to solve the set $\Xi$ of scenario subproblems.  Here, the parent node will allocate scenarios to the available children nodes dynamically (by keeping track of which ones are available and which
ones are busy solving another scenario subproblem). The children will communicate their solutions and cutting plane information to the parent node once they are done solving their subproblem.
The proposed graph abstraction allows us to simulate the effect of computing
and communication delays on the performance of the Benders scheme. This allows us, for instance, to simulate the behavior of the algorithm on a {\em hypothetical} parallel computer executing the algorithm that might
be subjected to a different number (and type) of CPUs or subjected to random computing loads from other jobs that increase computing latency in the CPUs.  Figure \ref{fig:benders_workflow_topology} depicts a hypothetical parallel computing
architecture with four computing nodes that execute the Benders algorithm. Here,
{\tt CPU 4} solves the master problem while {\tt CPU 1}, {\tt CPU 2}, and {\tt CPU 3} each solve a subproblem after receiving the scenario data attribute $\xi_n$ (where $n$ corresponds
to a CPU id) and the master solution attribute $\hat{x}$.

The simulation of the Benders algorithm  (Algorithm \ref{alg:synchronous_benders}) can be expressed in terms of nodes, tasks, and attributes following the setup provided in the Computing Graph \ref{alg:benders_workflow}.
The master node $m$ contains attributes defined in line \ref{line:benders_attributes_start} which include the solution to the master problem $\hat{x}$ and the scenario data $\{\xi_1...\xi_{N}\}$ which are communicated
to available subnodes $\mathcal{N}$.
The master node also contains tasks (line \ref{line:benders_master_tasks}) analogous to the functions in Algorithm \ref{alg:synchronous_benders} to execute the master problem ({$\tt run\_master$}) and to receive solutions from subnodes
(${\tt receive\_solution}$).  The details of each task and what attributes it updates can be found in Appendix \ref{app:benders}.
The {\tt solve\_master} task is triggered by the attribute $flag$ shown in line \ref{line:benders_run_master}.
The task solves problem \eqref{eq:benders_master}
using the master node attributes and \emph{updates} the attributes for the master solution $\hat{x}$ and scenarios $\{\xi_1...\xi_N\}$.   The updated
attributes trigger edge communication to the subnodes (lines \ref{line:edge_master_solution} and \ref{line:edge_scenario}).
The ${\tt receive\_solution}$ task is executed when any solution attribute $s_{n}$ is received from the connected
subnode $n$.  This task determines whether the $\tt run\_master$ task is ready to run again and either updates the attribute $flag$ (which triggers {\tt run\_master}) or updates the
corresponding scenario attribute $\xi_{n}$ with a new scenario to send back to subnode $n$.
Each subnode $n$ computes its task ${\tt solve\_subproblem}$ when it receives its scenario attribute $\xi$.
The ${\tt solve\_subproblem}$ task updates $s$, which triggers communication to the attribute $s_{n}$ on the master node $m$ (line \ref{line:sub_solution}).

\begin{footnotesize}
\makeatletter
\renewcommand*{\ALG@name}{Computing Graph}
\makeatother
\setcounter{algorithm}{0}

\begin{algorithm*}[h!]
\caption{(Synchronous Benders)}
\label{alg:benders_workflow}
\begin{algorithmic}[1]
    \State $\bm{Master \ Node} \ (m)$
    \State Attributes: $\bm{\mathcal{A}_{m}} := (\hat{x},S,C,flag,\{s_{1},...,s_{N}\},\{\xi_{1},...,\xi_{N}\})$ \label{line:benders_attributes_start}
    \State Tasks: $\bm{\mathcal{T}_{m}} := ({\tt run\_master},{\tt receive\_solution})$ \label{line:benders_master_tasks}
    \State \quad ${\tt run\_master}$: runs Task \ref{task:run_master}, triggered by: Updated($flag$) \label{line:benders_run_master}
    \State \quad ${\tt receive\_solution}$: runs Task \ref{task:receive_solution}, triggered by: Received($s_{n}$), \ $n \in \{1,...,N\}$
    \\\hrulefill
    \State $\bm{Sub \ Nodes} \ (n \in \mathcal{N})$
    \State Attributes: $\bm{\mathcal{A}_{n}} := (\hat{x},\xi,s)$
    \State Tasks: $\bm{\mathcal{T}_{n}} := ({\tt solve\_subproblem})$
    \State \quad ${\tt solve\_subproblem}$: runs Task \ref{task:run_sub_problem}, triggered by: Received($\xi$)
    \\\hrulefill
    \State $\bm{Edges \ \mathcal{E}_1} := \ \hat{x}_m \rightarrow \hat{x}_n,  \ n \in \{1,...,N\}$, send on: Updated($\hat{x}_m$) \label{line:edge_master_solution}
    \State $\bm{Edges \ \mathcal{E}_2} := \ \xi_{m,n} \rightarrow \xi_n, \ n \in \{1,...,N\}$, send on: Updated($\xi_{m,n}$) \label{line:edge_scenario}
    \State $\bm{Edges \ \mathcal{E}_3} := \ s_n  \rightarrow s_{m,n},  \ n \in \{1,...,N\}$, send on: Updated($s_n$) \label{line:sub_solution}
\end{algorithmic}
\end{algorithm*}
\end{footnotesize}

The implementation of the computing graph in ${\tt Plasmo.jl}$ for the case of three subnodes is shown in the code snippet of Figure \ref{fig:code_snippet_benders}.
Lines \ref{line:bendersworkflow1}-\ref{line:bendersworkflow2} create the graph, add the master node with its attributes and tasks (${\tt run\_master}$ and
${\tt receive\_solution}$).  Lines \ref{line:subnodes1}-\ref{line:subnodes2} add subnodes to the graph, each with a ${\tt solve\_subproblem}$ task that is
executed after receiving scenario data $\xi$ (line \ref{line:subnodetask}).
Finally, communication edges are created between attributes (lines \ref{line:connectb1}-\ref{line:connectb2}) and the graph
is executed (\ref{line:execute_benders}) until it terminates (i.e., a task calls ${\tt StopComputationGraph()}$ such as shown in Figure \ref{fig:code_snippet_benders_function}). It is also possible
to simulate to a pre-determined time by providing an argument to ${\tt execute!()}$.

Figure \ref{fig:code_snippet_benders_function} depicts how to define tasks in ${\tt Plasmo.jl}$.  One argument is typically needed: a reference to a node to retrieve and update attributes, but
a reference to the computing graph can be provided for access
to the graph clock or to terminate the computation.
For example, line \ref{line:cuts} retrieves the attribute value for the set of master problem cuts to
solve the master problem, and updates the master solution in line \ref{line:update_solution}.
Also note that the ${\tt run\_master}$ and ${\tt solve\_subproblem}$ tasks are given default compute times as the \emph{true compute time} of their execution (${\tt compute\_time = :walltime}$).
We set a communication delay of 0.005 seconds from the master
node to the subnodes (but it is also possible to make delay time a function of the attributes communicated or to experiment with different delays to evaluate the effect of
communication overhead).

\begin{figure}[h!]
\centering
\begin{scriptsize}
\lstset{language=Julia}
\begin{lstlisting}[escapeinside={(*}{*)},escapechar = |]
# function implementing master task.
function run_master(graph::ComputingGraph,node::ComputeNode)
    C = getattribute(node,:C)  |\label{line:cuts}|
    S = getattribute(node,:S)
    solution = solve_master_problem(C)
    updateattribute(node[:solution],solution)   |\label{line:update_solution}|

    #Convergence check.
    lower_bound = solution.objval           #compute the lower bound given the master solution
    upper_bound = compute_upper_bound(S)    #compute the upper bound from sub-problem solutions
    if converged(lower_bound,upper_bound)
        StopComputingGraph(graph)         #stop the computing graph if converged
    end
    |$\Xi$| = getattribute(node,:|$\Xi$|)      #set of scenarios to evaluate
    #Start sending out scenarios.  Each scenario update will trigger communication a subnode.
    |$\xi$| = getattribute(node,:|$\xi$|) #retrieve array of outgoing scenarios
    for i = 1:N
        updateattribute(|$\xi$|[i],|$\Xi$[i]|)
    end
end
\end{lstlisting}
\end{scriptsize} %
\caption{${\tt Plasmo.jl}$ snippet demonstrating a template task for a computing graph.}
\label{fig:code_snippet_benders_function}
\end{figure}

\begin{figure}[h!]
\centering
\begin{scriptsize}
\lstset{language=Julia}
\begin{lstlisting}[escapeinside={(*}{*)},escapechar = |]
    #Create a computing graph
    graph = ComputingGraph()            |\label{line:bendersworkflow1}|

    N = 3  #number of subnodes
    #Add the master node (m)
    m = addnode!(graph)
    @attributes(m,x,C,flag,|$\xi$|[1:N],s[1:N])
    @nodetask(graph,m,run_master,compute_time = :walltime, triggered_by=Updated(flag))
    @nodetask(graph,m,receive_solution[i = 1:N], compute_time = 0, triggered_by=Received(s[i]))
    #Provide an initial signal for the computation-graph to evaluate
    schedule_trigger(graph,run_master, time = 0)  |\label{line:bendersworkflow2}|

    for i = 1:N                |\label{line:subnodes1}|
        #Add subnodes (n = 1:N) to solve sub-problems
        subnode = addnode!(graph)
        @attributes(subnode,x,|$\xi$|,s)
        @nodetask(graph,subnode,solve_subproblem, compute_time = :walltime,triggered_by=Received(|$\xi$|))  |\label{line:subnodetask}|

        #Connect attributes between master and subnodes
        @connect(graph,m[:x] => n[:x],send_on=Updated(m[:x]))   |\label{line:connectb1}|
        @connect(graph,m[:|$\xi$|][i] => n[:|$\xi$|],send_on=Updated(m[:|$\xi$|][i]),delay=0.005)
        @connect(graph,n[:s] => m[:s][i],send_on=Updated(n[:s])) |\label{line:connectb2}|
    end                                     |\label{line:subnodes2}|
    #Execute the computing graph
    execute!(graph) |\label{line:execute_benders}|
\end{lstlisting}
\end{scriptsize} %
\caption{Plasmo.jl snippet for simulating Benders decomposition.}
\label{fig:code_snippet_benders}
\end{figure}

\begin{figure}[h]
    \centering
    \includegraphics[scale=0.7]{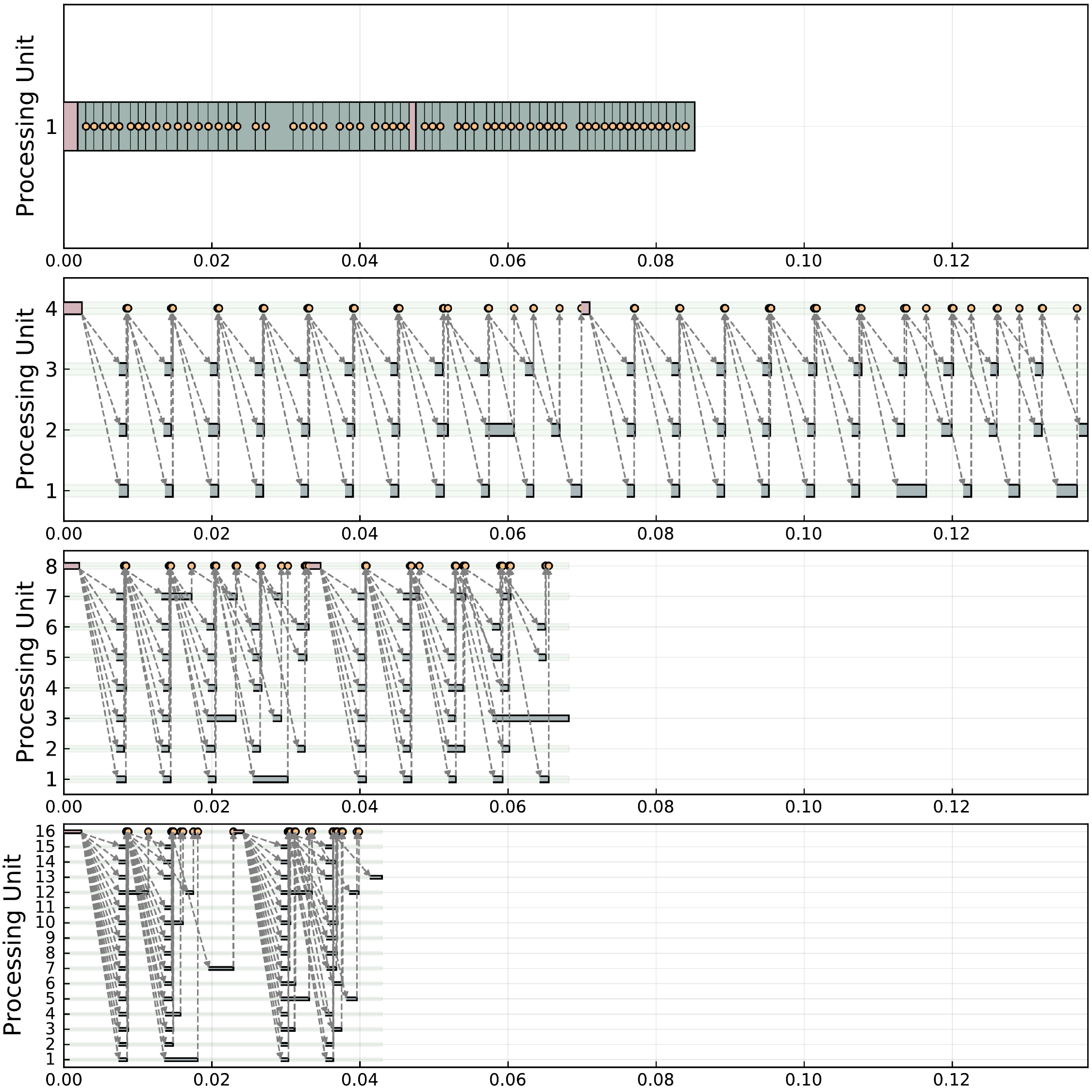}
    \caption{Simulation of Benders decomposition using different numbers of CPUs (top panel shows one CPU while bottom panel shows sixteen CPUs).  Red tasks correspond to the master problem execution time,
    grey tasks represent subproblem execution time and orange dots represent the {\tt receive\_solution} task which is simulated with zero computing time.}
    \label{fig:benders_workflow_sim}
\end{figure}

Figure \ref{fig:benders_workflow_sim} summarizes the simulation results of the Benders algorithm as we increase the number of CPUs available in the computing architecture (we consider cases with $N=1,4,8,$ and $16$ CPUs).
We can see that, with a communication delay of 0.005 seconds from the master to the subnodes, using one CPU has a shorter total solution time than using four CPU nodes (due to the
communication overhead).  Executing the algorithm on 8 and 16 CPU nodes, however,
results in algorithm speed up (reduction in computing latency overcomes communication latency).
This illustrates how the proposed framework can help predict trade-offs of computing and communication latency.
For instance, our results predict that the proposed Benders scheme only benefits from parallelization when the number of CPUs is sufficiently large. We highlight that the parallel architectures evaluated are {\em hypothetical} (the actual simulation of the algorithm was executed on a single CPU). In other words, ${\tt Plasmo.jl}$ simulates the behavior of the Benders algorithm on a virtual computing environment.

\FloatBarrier

\subsection{Model Predictive Control Architectures}

\begin{figure}[h]
    \centering
    \includegraphics[scale=0.5]{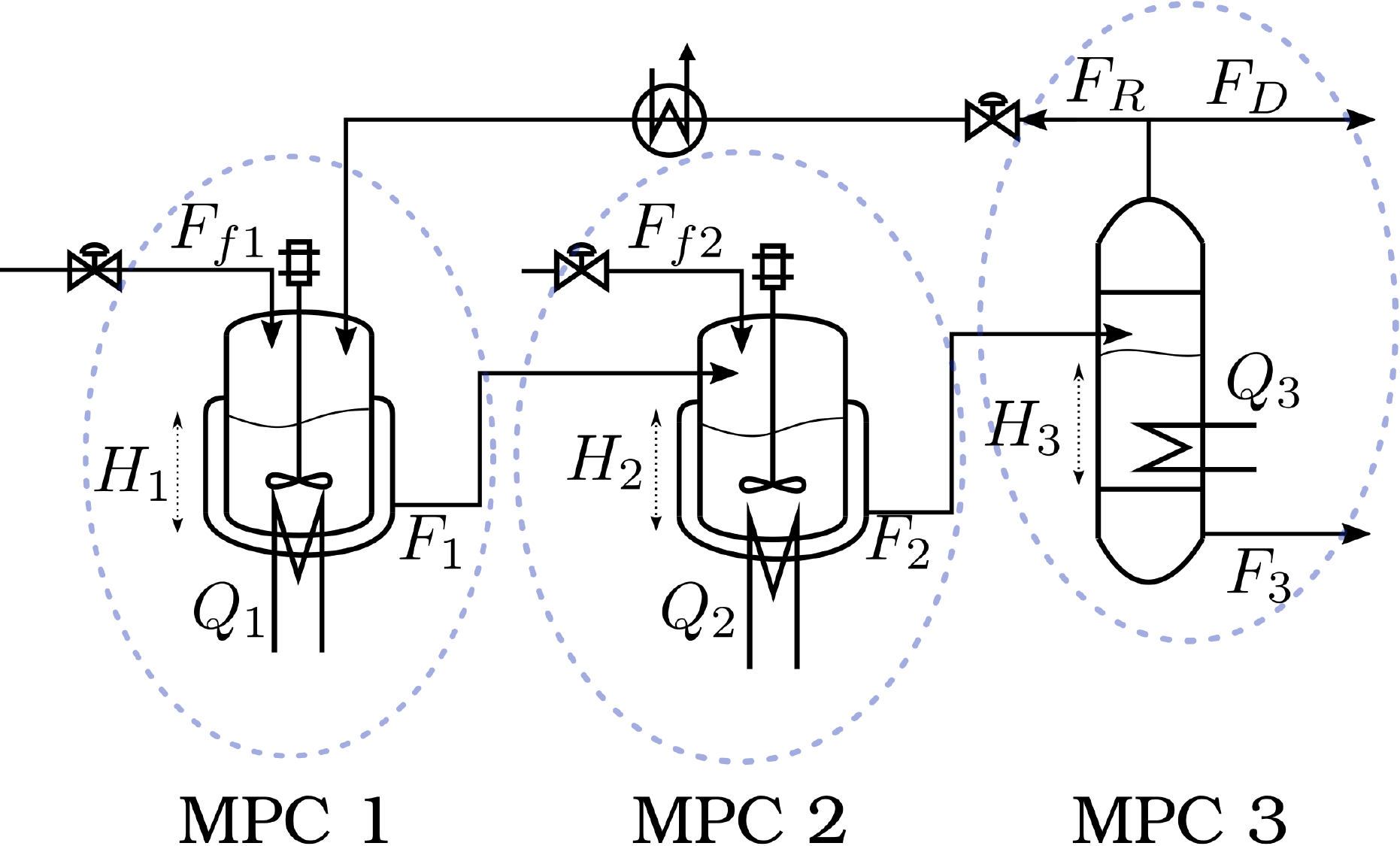}
    \caption{Reactor separator process and partitioning into MPC controllers.}
    \label{fig:reactor_diagram}
\end{figure}

This case study demonstrates how a computing graph can be used to simulate the behavior of distributed control
architectures when considering computation and communication delays and to simulate how such delays impact the actual behavior of a physical system. Specifically, we consider a reactor-separator system (see Figure \ref{fig:reactor_diagram})
from \cite{Stewart2010}, which is a standard application for evaluating distributed model predictive control (MPC) algorithms.
The system consists of two reactors in series where the reaction $A \ \rightarrow \ B \ \rightarrow \ C$ takes place and a separator which produces
a product stream and a recycle.  The system is described by twelve states: the weight fractions of A and B in each unit, the unit heights, and the unit temperatures.  The manipulated
inputs are the flow rates $F_{f1}$, $F_{f2}$, $F_1$, $F_2$, $F_3$, $F_R$, and heat exchange rates $Q_1$, $Q_2$, and $Q_3$.  The details of the system are given in Appendix \ref{app:reactor}.

\begin{figure}[h]
\centering
\begin{subfigure}[t]{0.33\textwidth}
    \centering
    \includegraphics[scale=0.40]{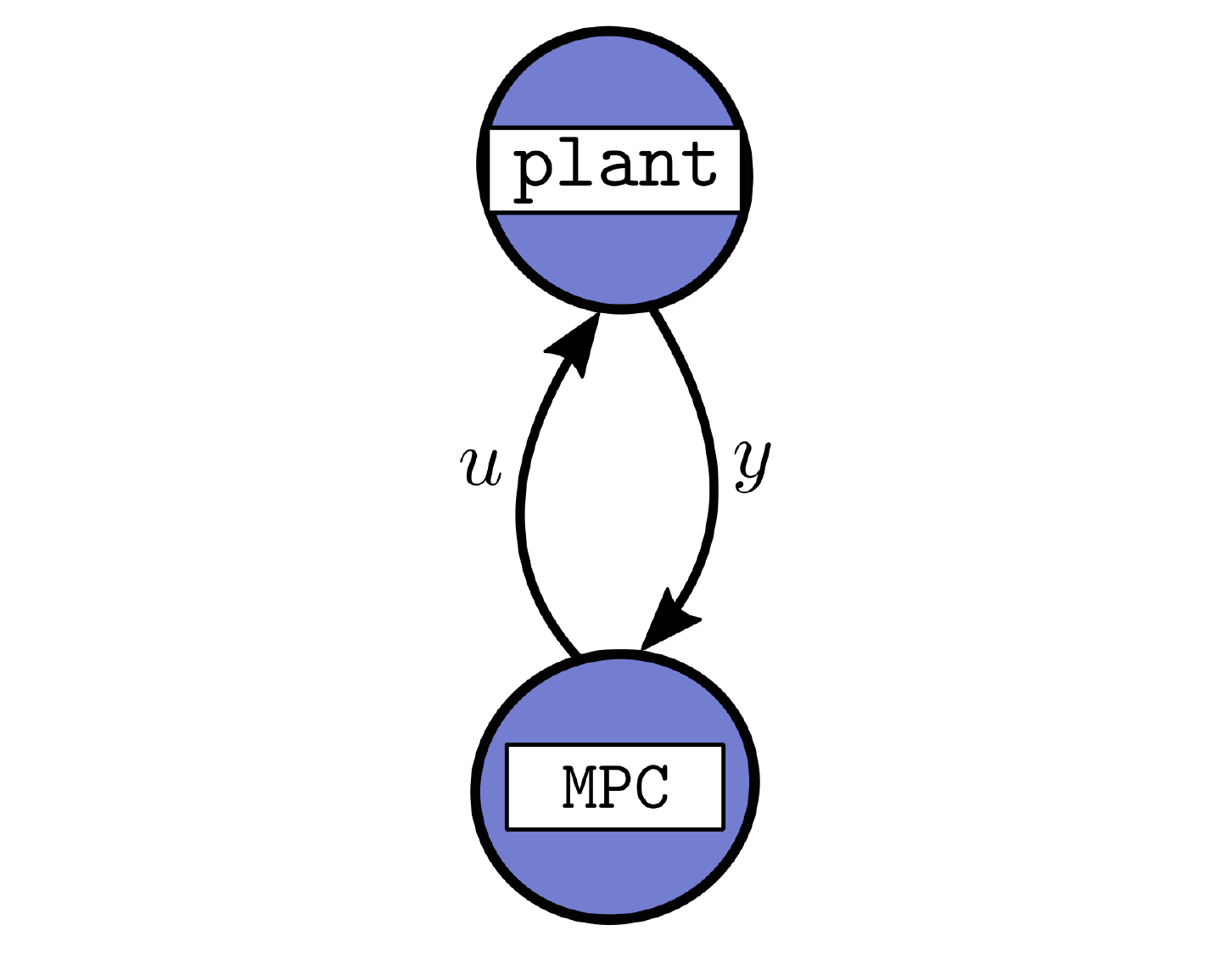}
    \caption{Centralized}
    \label{fig:centralized_topology}
\end{subfigure}
\begin{subfigure}[t]{0.33\textwidth}
    \centering
    \includegraphics[scale=0.40]{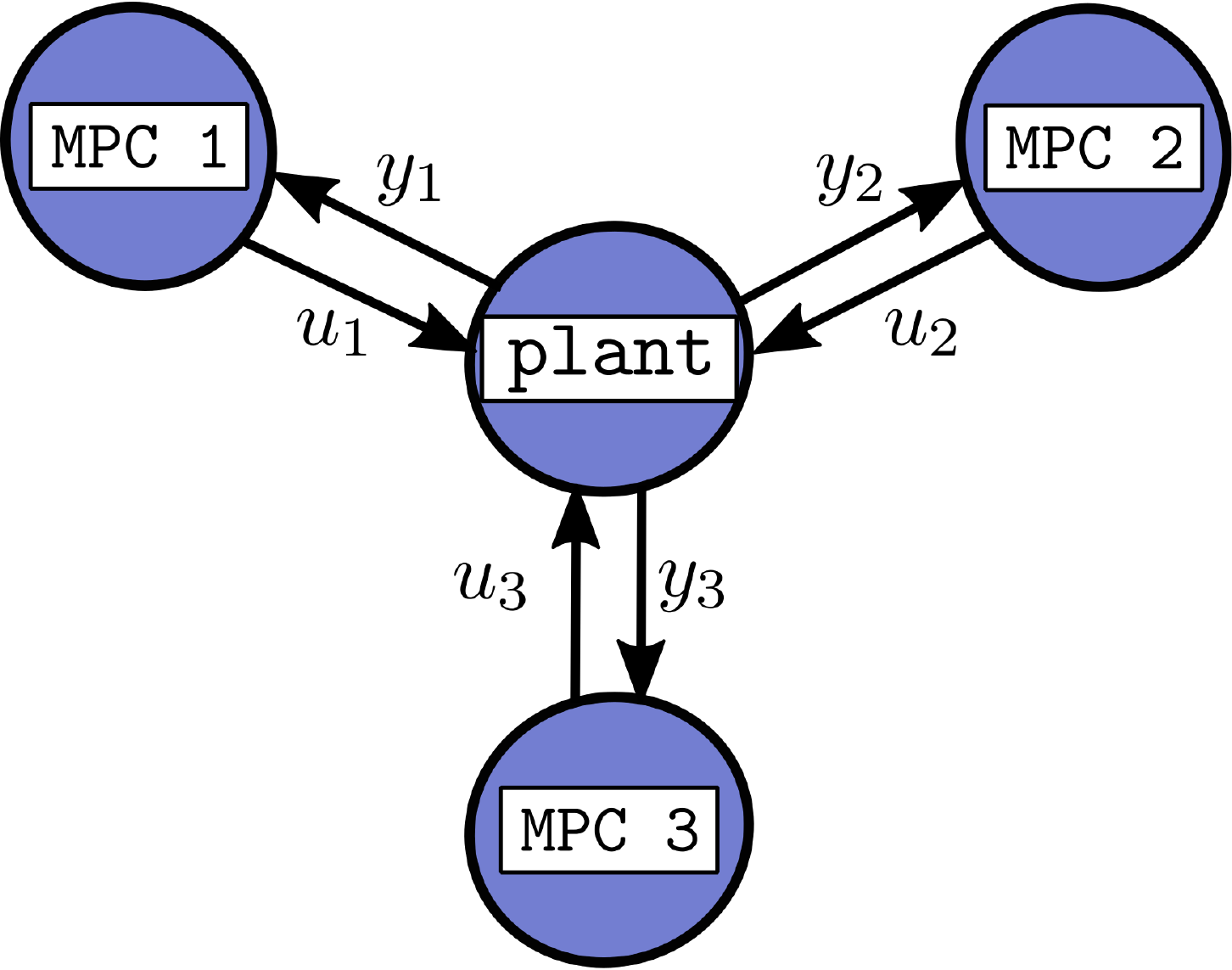}
    \caption{Decentralized}
    \label{fig:decentralized_topology}
\end{subfigure}
\\
\vspace{0.2in}
\begin{subfigure}[t]{0.33\textwidth}
    \centering
    \includegraphics[scale=0.40]{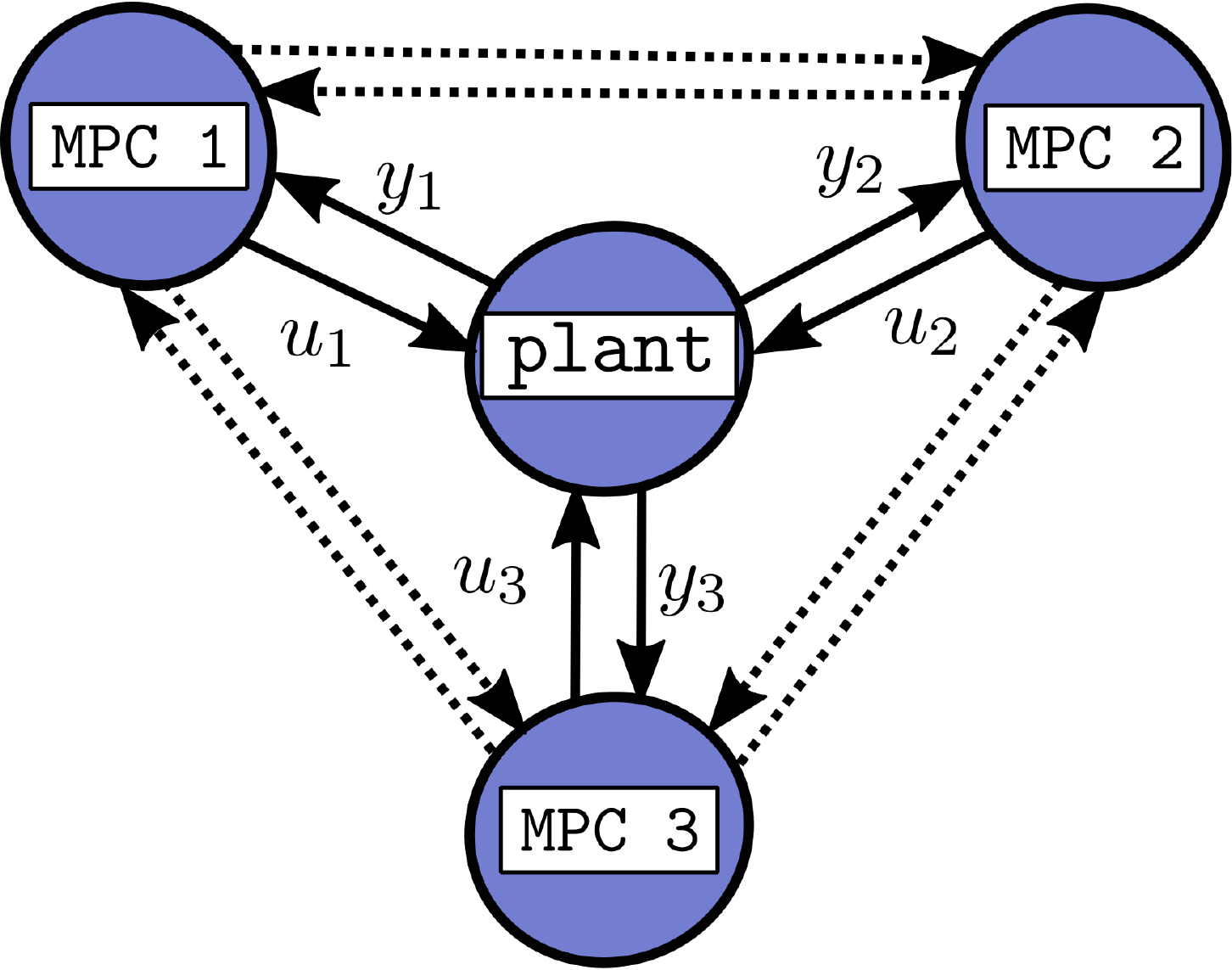}
    \caption{Cooperative}
    \label{fig:cooperative_topology}
\end{subfigure}
\caption{Simulated MPC architectures:  centralized (top left),  decentralized (top right) and cooperative (bottom).}
\end{figure}

The behavior of the physical system (the plant) shown in Figure \ref{fig:reactor_diagram} is simulated under three different MPC architectures.
We consider a \emph{centralized} MPC architecture (Figure \ref{fig:centralized_topology}) wherein every output is sent to a central MPC controller which computes all control actions for the plant.
We also consider \emph{decentralized} control architectures consisting of three MPC controllers (one for each unit) and simulate their behavior when they do not
communicate (Figure \ref{fig:decentralized_topology}) and when they \emph{cooperate} by communicating state and intended control actions
(Figure \ref{fig:cooperative_topology}). Complex performance, computation, and communication trade offs arise under the three MPC architectures studied. In particular, the centralized scheme achieves best performance when the computing and communication delays are short (which might not be achievable in large systems). On the other hand, the performance of decentralized schemes might be lower but computing delays are expected to be shorter as well. Analyzing such trade offs is facilitated by the computing graph, since this captures communication and computing times while simultaneously advancing the plant simulation. The proposed framework also captures \emph{asynchronous} behavior of the decentralized and cooperative schemes.
In particular, the controllers might inject their control actions as soon as they complete their computing task (as opposed to waiting when all of them are done).

Computing Graph \ref{alg:cooperative_workflow} specifies the setup for the cooperative MPC algorithm.  The plant node $pl$ contains the task ${\tt run\_plant}$ which
advances the state of the system from the current clock time to the time of the next action signal.  The plant node includes attributes $u_1$, $u_2$, and $u_3$ (which are the
control actions received from the MPC controllers) and $x$ (which are the plant states).  The plant node
is self-triggered by updating $x_{pl}$, which allows the task to run continuously (i.e., the task ${\tt run\_plant}$ constantly updates the attribute $x_{pl}$).  The plant state is communicated
to each MPC node at a {\em constant sampling time}.  The edges that connect the plant state attribute $x_pl$ to each MPC node will trigger when sending its source attribute (with a given waiting time $\theta_{wait}$).
The MPC nodes $n_{mpc,1}$, $n_{mpc,2}$, and $n_{mpc,3}$ each execute their task ${\tt control\_action}$ which computes a control action by solving an optimization problem using attributes from the rest of the
MPC controllers and from the plant.  If they have completed enough
coordination iterations, they each update their attribute $u_{inject}$, which triggers communication into the plant.
Otherwise, they each update their control trajectory $u_{p_n}$ and
exchange their attributes. This triggers the MPC node task ${\tt receive\_policy}$, which manages the MPC trajectory exchange.
Additionally, we specify logic to handle the case when the attributes \{$u_{p_i}, (i \neq n)$\} are received while a node $n$ is busy executing a task.
When this occurs, the triggered task is queued and re-triggered when the current task being executed is completed.

\makeatletter
\renewcommand*{\ALG@name}{Computing Graph}
\makeatother

\begin{algorithm*}[h!]
\caption{Cooperative Control}
\label{alg:cooperative_workflow}
\begin{algorithmic}[1]
    \State $\bm{Plant \ Node} \ (pl)$
    \State Attributes: $\bm{\mathcal{A}_{pl}} := (u_{1},u_{2},u_{3},x)$
    \State Tasks: $\bm{\mathcal{T}_{pl}} := ({\tt run\_plant})$
    \State $\quad {\tt run\_plant}$: Task \ref{task:run_plant}, triggered by: Updated($x$)
    \\\hrulefill
    \State $\bm{MPC \ Nodes} \ (n = \{1,2,3\})$
    \State Attributes: $\bm{\mathcal{A}_{n}} := (u_{inject},y,u_{p_1},u_{p_2},u_{p_3},iter,flag)$
    \State Tasks: $\bm{\mathcal{T}_{n}} := ({\tt control\_action},{\tt receive\_policy})$
    \State \quad ${\tt control\_action}$: Task \ref{task:control_action}, triggered by: Received($y$) or Updated($flag$)
    \State \quad ${\tt receive\_policy}$: Task \ref{task:receive_policy}, triggered by: Received($u_{p_i}), \ i \neq n$
    \\\hrulefill
    \State $\bm{Edges \ \mathcal{E}_1} := \ x_{pl} \rightarrow y_n, \ n = \{1,2,3\}$, send on: Sent($x_{pl}$), wait: $\theta_{wait}$
    \State $\bm{Edges \ \mathcal{E}_2} := \ u_{n,inject} \rightarrow u_{pl,n}, \ n = \{1,2,3\}$, send on: Updated($u_{n,inject}$)
    \State $\bm{Edges \ \mathcal{E}_3} := \ u_{n,p_n} \rightarrow u_{i,p_n}, \ n = \{ 1,2,3 \}, i \neq n$, send on: Updated($u_{n,p_n}$)
\end{algorithmic}
\end{algorithm*}

Figure \ref{fig:cooperative_code_snippet} demonstrates how the cooperative MPC algorithm (Computing Graph \ref{alg:cooperative_workflow}) is implemented in ${\tt Plasmo.jl}$.
A node is created for the plant and for each MPC controller.  Tasks and attributes are added to each node and attributes for sensor measurements and control inputs are connected between nodes.
The execution behaviors are once again modified using keyword arguments.  Lines \ref{line:coop_start}
to \ref{line:coop_end} create the graph, add a node for the plant simulation, and a node for each MPC controller.  Task execution is
triggered by receiving and updating attributes.  Each MPC controller computes its control action when it samples (i.e., receives the output measurement attribute $y$) and starts
performing iterations.  The measurement attributes $y$ are communicated to each controller with a delay of
30 seconds at a sample period of 60 seconds (by specifying ${\tt send\_wait = 60}$), starting at 5 seconds (line \ref{line:connectcoop1}).  The controller inputs are connected
to the corresponding plant attributes and are communicated when they are updated (line \ref{line:connectcoop2}).  Finally, lines \ref{line:connectcoop3}-\ref{line:connectcoop4}
connect the MPC controllers to each other (create communication edges) so that they exchange information when they update their control actions.

\begin{figure}[h!]
\centering
\begin{scriptsize}
\lstset{language=Julia}
\begin{lstlisting}[escapeinside={(*}{*)},escapechar = |]
    #Create a computing graph
    graph = ComputationGraph()                                |\label{line:coop_start}|

    #Create plant simulation node
    pl = addnode!(graph)
    @attributes(pl,x,u1,u2,u3)
    @nodetask(graph,pl,run_plant,triggered_by=Updated(x))
    schedule_trigger(graph,run_plant,time = 0) #initialize computing graph

    #Reactor 1 MPC
    n1 = addnode!(graph)
    @attributes(n1,u_inject,y,u_p1,u_p2,u_p3,iter,flag)
    @nodetask(graph,n1,control_action_r1,triggered_by=[Received(y),Updated(flag)])
    @nodetask(graph,n1,receive_policy,triggered_by=Received(u_p2,u_p3),trigger_during_busy=:queue_task)

    #Reactor 2 MPC
    n2 = addnode!(graph)
    @attributes(n2,u_inject,y,u_p1,u_p2,u_p3,iter,flag)
    @nodetask(graph,n2,control_action_r2,triggered_by=[Received(y),Updated(flag)])
    @nodetask(graph,n2,receive_policy,triggered_by=Received(u_p1,u_p3),trigger_during_busy=:queue_task)

    #Separator MPC
    n3 = addnode!(graph)
    @attributes(n3,u_inject,y,u_p1,u_p2,u_p3,iter,flag)
    @nodetask(graph,n3,control_action_sep,triggered_by = [Received(y),Updated(flag)])
    @nodetask(graph,n3,receive_policy,triggered_by=Received(u_p1,u_p2),trigger_during_busy=:queue_task)   |\label{line:coop_end}|

    #Connect plant to MPC controllers
    @connect(graph,pl[:x]=>[n1[:y],n2[:y],n3[:y]],delay=30,send_on=Sent(pl[:x]),send_wait=60,start=5)   |\label{line:connectcoop1}|
    #Connect MPC controllers to plant
    @connect(graph,[n1[:u_inject],n2[:u_inject],n3[:u_inject]]=>[pl[:u1],pl[:u2],pl[:u3]],delay=30,send_on=Updated([n1[:u_inject],n2[:u_inject],n3[:u_inject]])) |\label{line:connectcoop2}|

    #Connect MPC controllers to perform cooperation
    @connect(graph,n1[:u_p1]=>[n2[:u_p1],n3[:u_p1]],send_on=Updated(n1[:u_p1]))     |\label{line:connectcoop3}|
    @connect(graph,n2[:u_p2]=>[n1[:u_p2],n3[:u_p2]],send_on=Updated(n2[:u_p2]))
    @connect(graph,n3[:u_p3]=>[n1[:u_p3],n2[:u_p3]],send_on=Updated(n3[:u_p3]))     |\label{line:connectcoop4}|

    #Execute computing graph
    execute!(graph,5000)
\end{lstlisting}
\end{scriptsize} %
\caption{Plasmo.jl code snippet for simulating the cooperative MPC architecture.}
\label{fig:cooperative_code_snippet}
\end{figure}

\begin{figure}[h!]
\centering
\begin{subfigure}[b]{0.48\textwidth}
    \centering
    \includegraphics[scale=0.6]{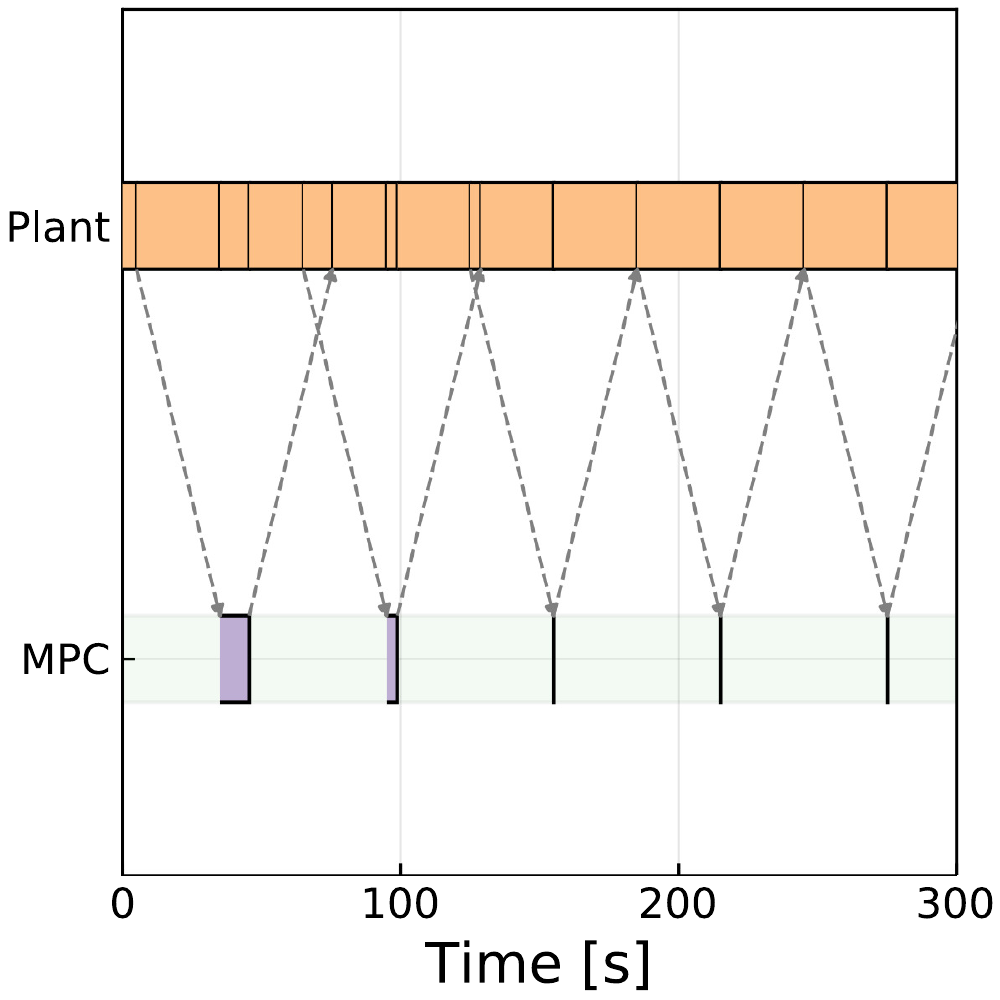}
    \caption{Centralized communication}
    \label{fig:centralized_pattern}
\end{subfigure}
\begin{subfigure}[b]{0.48\textwidth}
    \centering
    \includegraphics[scale=0.6]{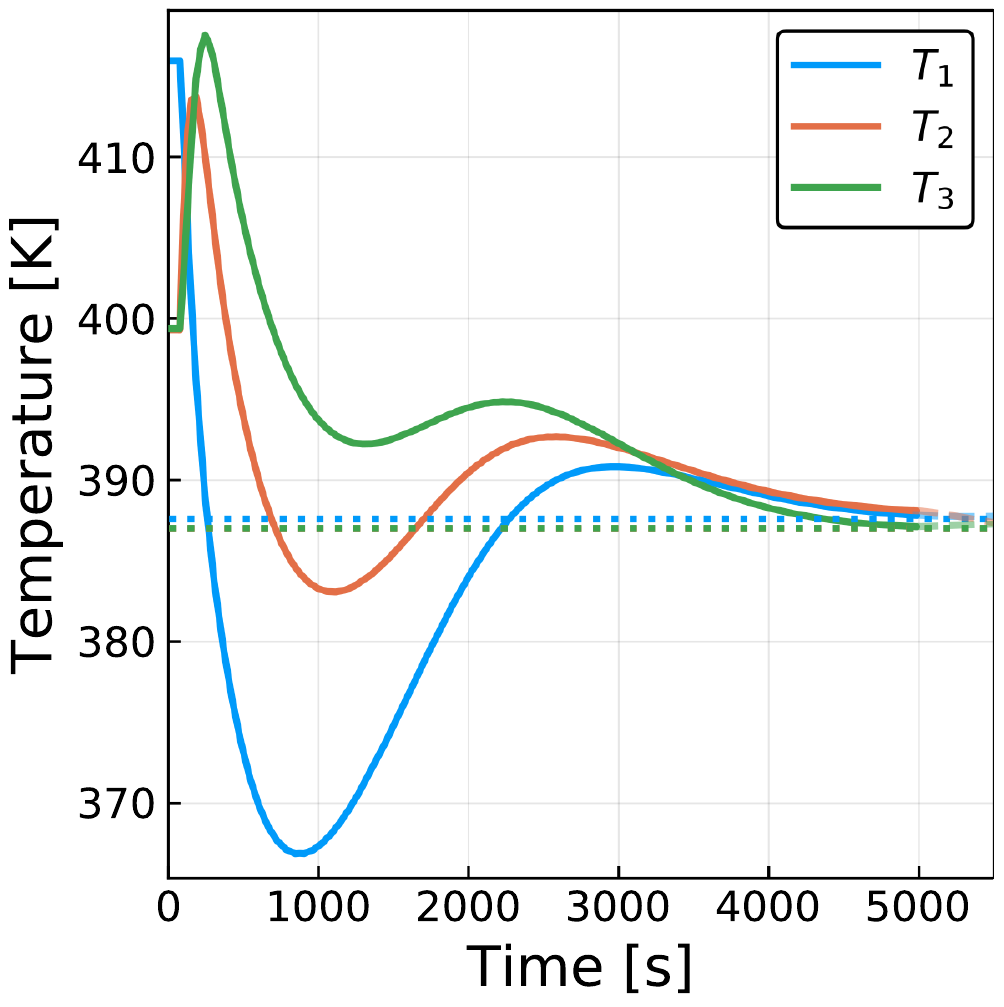}
    \caption{Centralized state profile}
    \label{fig:centralized_temp}
\end{subfigure}
\\
\begin{subfigure}[b]{0.48\textwidth}
    \centering
    \includegraphics[scale=0.6]{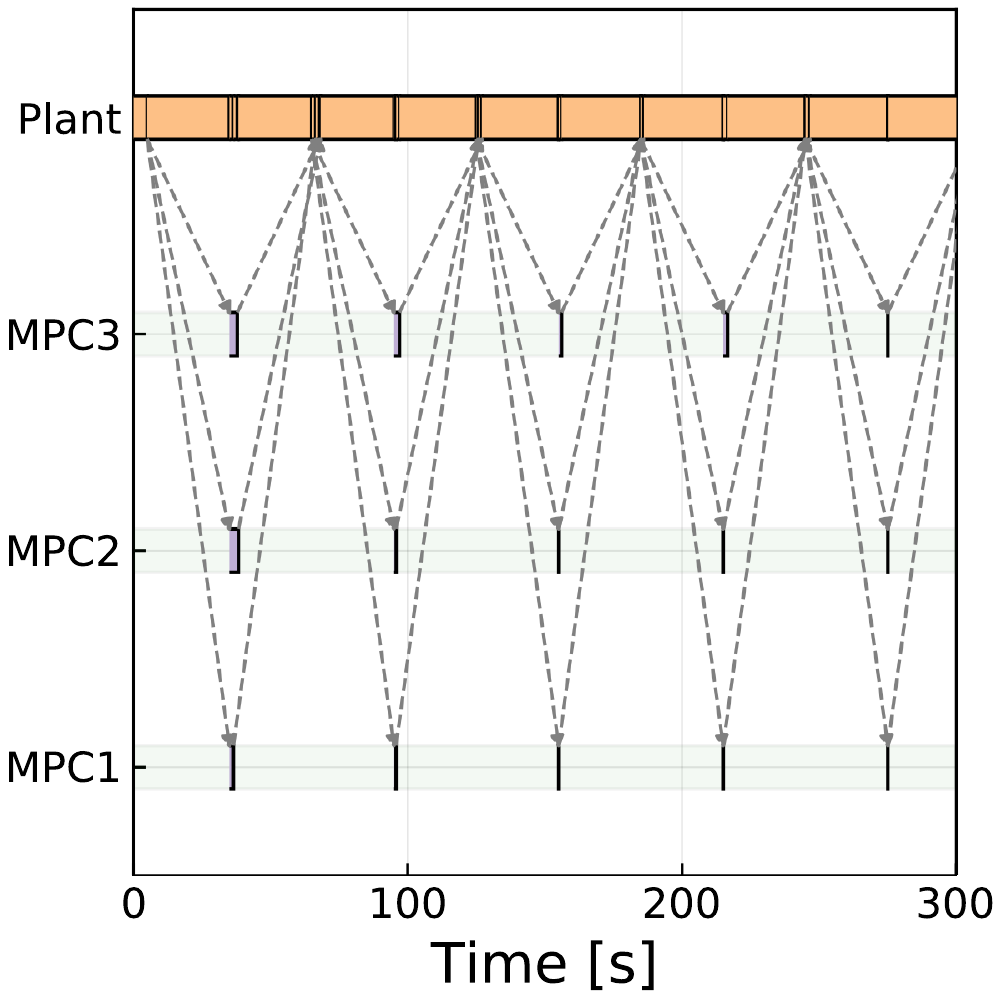}
    \caption{Decentralized communication}
    \label{fig:decentralized_pattern}
\end{subfigure}
\begin{subfigure}[b]{0.48\textwidth}
    \centering
    \includegraphics[scale=0.6]{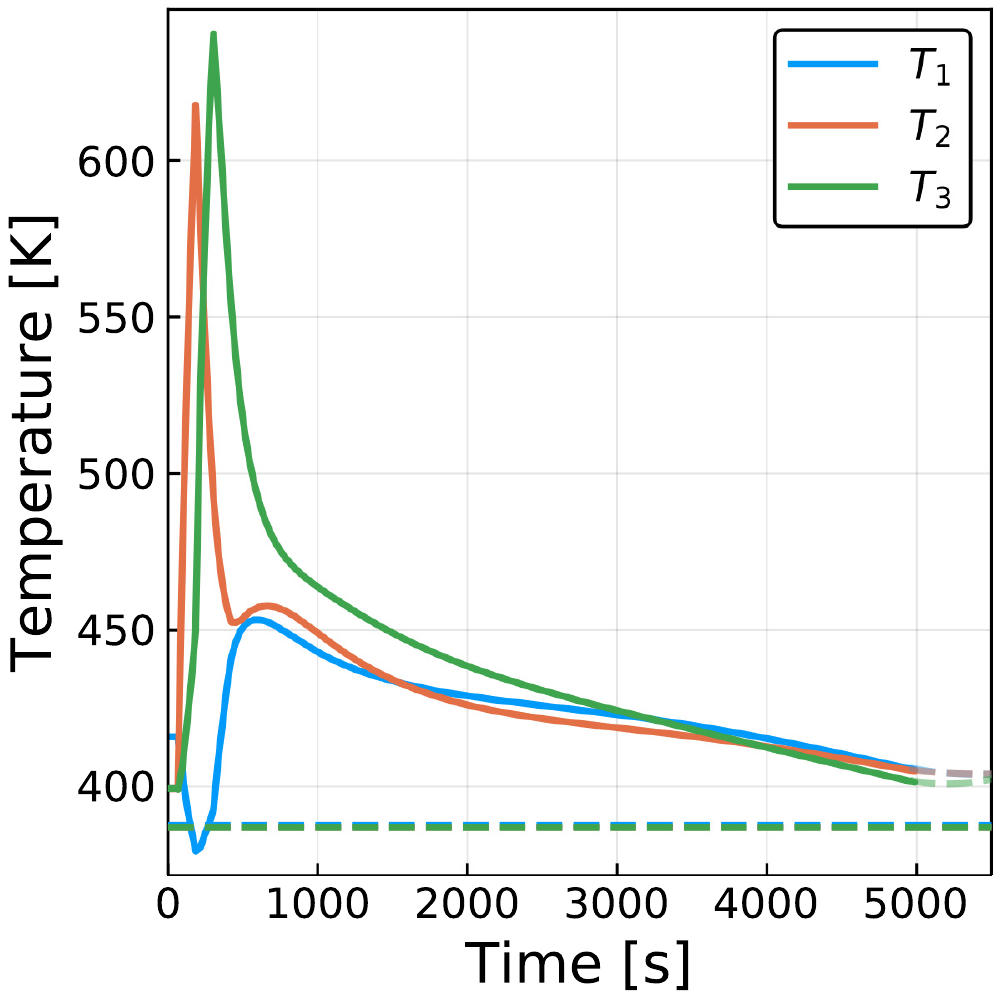}
    \caption{Decentralized state profile}
    \label{fig:decentralized_temp}
\end{subfigure}
\\
\begin{subfigure}[b]{0.48\textwidth}
    \centering
    \includegraphics[scale=0.6]{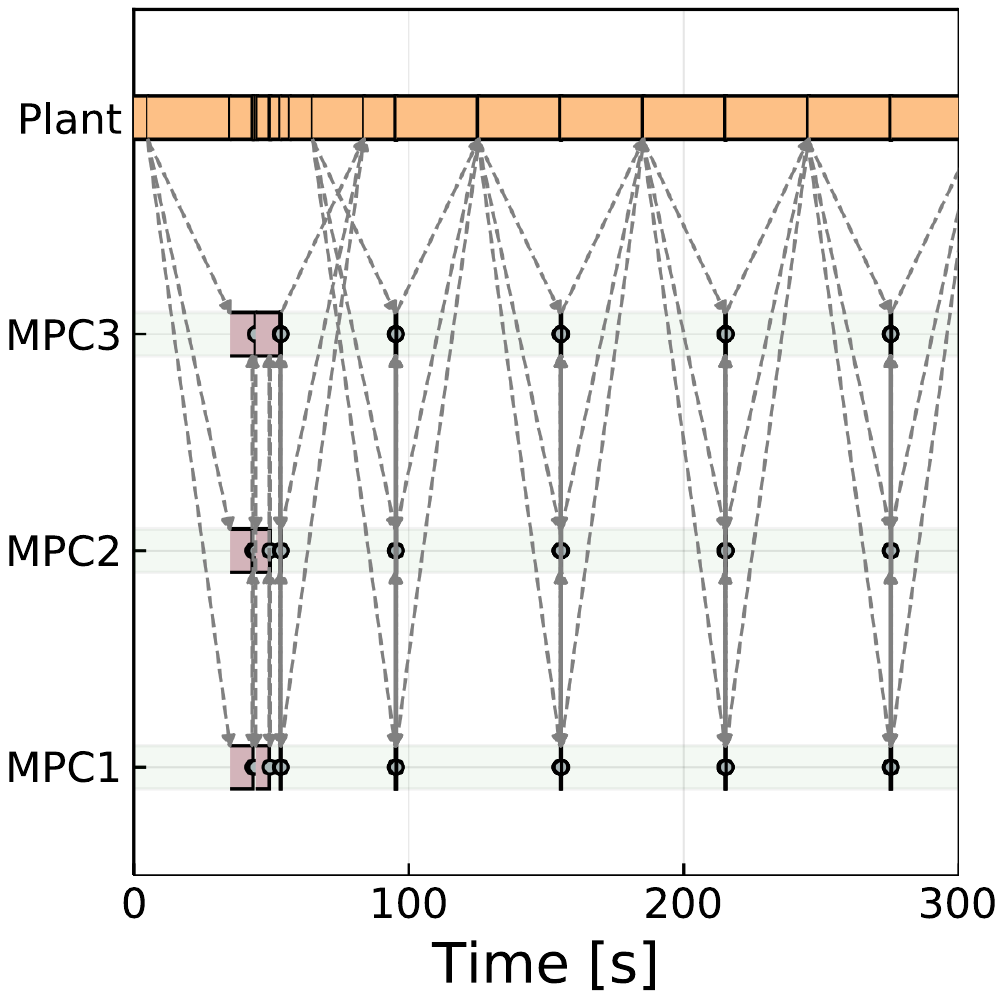}
    \caption{Cooperative communication}
    \label{fig:cooperative_pattern}
\end{subfigure}
\begin{subfigure}[b]{0.48\textwidth}
    \centering
    \includegraphics[scale=0.6]{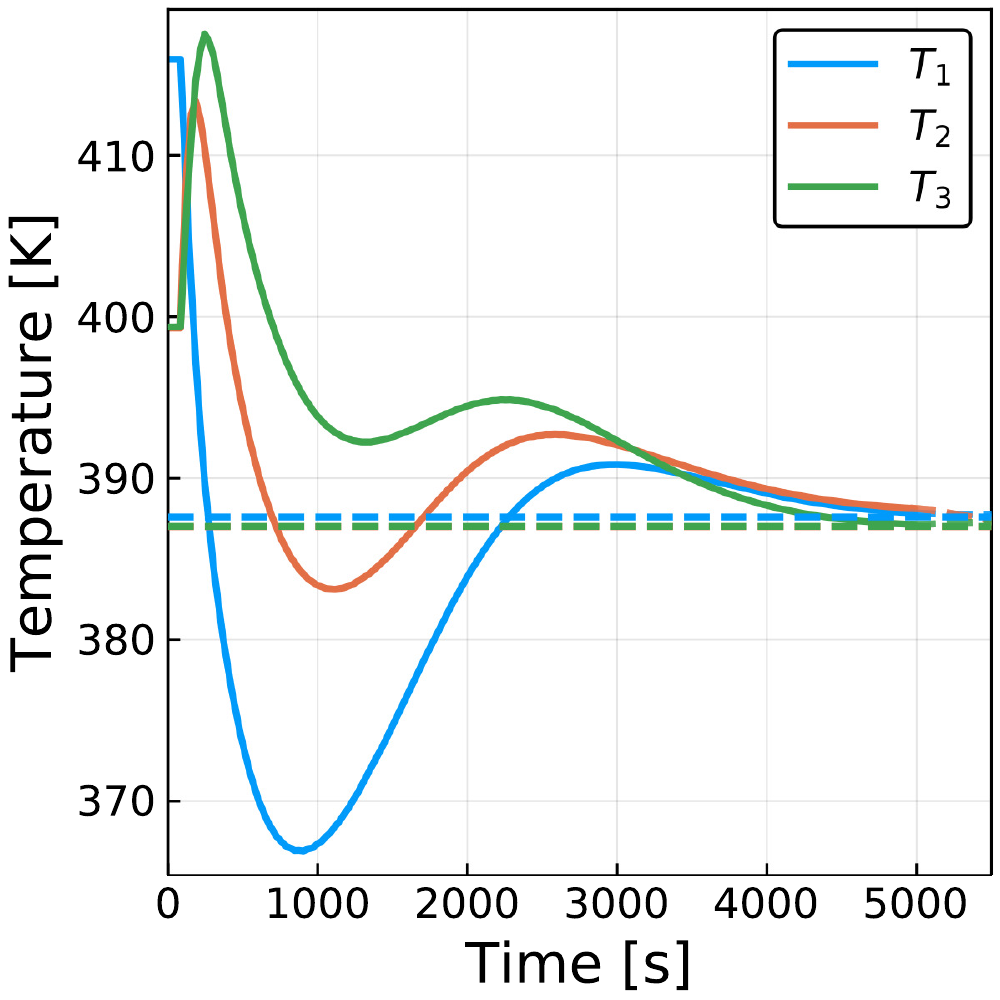}
    \caption{Cooperative state profile}
    \label{fig:cooperative_temp}
\end{subfigure}
\caption{Simulation results for MPC architectures.  Centralized MPC converges to the set-point despite the computing delays (top panels). Decentralized
MPC does not converge to the set-point (middle panels).  Cooperate MPC exhibits communication complexity but converges to the centralized MPC solution (bottom panels).}
\label{fig:distributed_results}
\end{figure}

Figure \ref{fig:distributed_results} presents the simulation results for each MPC algorithm.  The centralized MPC communication pattern (\ref{fig:centralized_pattern}) shows the
communication delays between the plant and the controller (grey arrows),
the time required to compute the control action (the purple bar), and highlights how the plant
state advances continuously while computation and communication tasks execute.  Despite the delays enforced for the controller, centralized
MPC is able to drive the state to the set-point (Figure \ref{fig:centralized_temp}).  Decentralized MPC does not require communication
between controllers and computing times are decreased (Figure \ref{fig:decentralized_pattern}) but we observe
that the set-point cannot be reached (Figure \ref{fig:decentralized_temp}).  This is because this approach does not adequately capture multi-variable interactions.
Finally, cooperative MPC has a more sophisticated communication strategy (Figure \ref{fig:cooperative_pattern})
but we observe that this helps mimic the performance of centralized MPC (Figure \ref{fig:cooperative_temp}).

\FloatBarrier

\section{Conclusions and Future Work}
We have presented graph-based abstractions that facilitate modeling complex cyber-physical systems.  An algebraic model graph abstraction is used to facilitate the construction,
solution, and analysis of complex optimization models while a computing graph abstraction is used to facilitate the construction, solution, and analysis of complex decision-making algorithms that involve diverse computing tasks and communication protocols. We implemented both abstractions in the Julia package ${\tt Plasmo.jl}$ and provided
case studies that demonstrate the capabilites of the abstractions.

Future work will focus on developing and implementing more sophisticated graph analysis techniques for the model graph that facilitate decomposition strategies and visualization.
We are also interested in extending the computing graph capabilties to allow for co-simulation \cite{Gomes2018},  which would facilitate
decomposition of the simulation framework. This is necessary to enable scalability to systems with many computing nodes such as swarm-based systems (which contain thousands of nodes) \cite{Lopes2016,Scott2017} or control systems with adaptive communication \cite{Bodson2011}. We are also interested in using the developed capabilities to simulate more sophisticated hierarchical control architectures that involve combinations of centralized/decentralized subsystems \cite{shinzavala}.

\section*{Acknowledgements}
This material is based on work supported by the U.S. Department of Energy (DOE), Office of Science, under Contract No. DE-AC02-06CH11357 as well as the DOE Office of Electricity Delivery and Energy Reliability’s Advanced Grid Research and Development program (AGR\&D).  This work was also partially supported by the U.S. Department of Energy grant DE-SC0014114. 

\newpage

\appendix

\section{Gas Pipeline Study Model}

We model gas pipeline networks as algebraic graphs $\mathcal{G}(\mathcal{N},\mathcal{E})$ wherein the nodes $\mathcal{N}$ consist of models for
the links (pipelines) in the network $\mathcal{L}_{gas} \subseteq \mathcal{N}$ and the gas network junctions $\mathcal{N}_{gas} \subseteq \mathcal{N}$.
By defining component models for each node in the graph, we express the physical system as a collection of models connected by algebraic link constraints.

\subsection{Gas Pipeline Component  Model}

We assume isothermal flow through a horizontal pipeline segment $\ell \in \mathcal{L}_{gas}$. The isothermal mass and momentum conservation equations
are presented in \cite{Osiadacz1984} and take the form:
\begin{subequations}\label{eq:isothermaleuler}
    \begin{align}
        &\frac{\partial \rho_\ell(t,x)}{\partial t} + \frac{\partial (\rho_\ell(t,x) v_\ell(t,x))}{\partial x} = 0\\
        &\frac{\partial (\rho_\ell(t,x) v_\ell(t,x))}{\partial t} + \frac{\partial(\rho_\ell(t,x) v_\ell(t,x)^2+ p_\ell(t,x))}{\partial x} =
        -\frac{\lambda_\ell}{2D}\rho_\ell(t,x) v_\ell(t,x)\mods{v_\ell(t,x)},\label{eq:isothermalmomentum}
     \end{align}
\end{subequations}
where notation and units are defined in Table \ref{table:variables}.
The euler equations can be approximated by dropping the the momentum term
$\partial_x(\rho v^2)$ which has a negligible effect on dynamics, and the states $\rho$ and $v$ can be converted into mass flow rate $f$ and pressure $p$ to produce
the following formulation
\begin{subequations}\label{eq:weymouth}
    \begin{align}
            \frac{\partial p_\ell(t,x)}{\partial t} + c_{1,\ell}\frac{\partial f_\ell(t,x)}{\partial x} &= 0\label{eq:weymouthmass}\\
            \frac{\partial f_\ell(t,x)}{\partial t} + c_{2,\ell}\frac{\partial p_\ell(t,x)}{\partial x} + c_{3,\ell}\frac{f_\ell(t,x)}{p_\ell(t,x)} \mods{f_\ell(t,x)}&= 0,\label{eq:weymouthmomentum}
    \end{align}
\end{subequations}
where  $c_{1,\ell}$,$c_{2,\ell}$ and $c_{3,\ell}$ are defined constants in Table \ref{table:parameters}.

The flow and pressure variables can also be used to define boundary conditions for the link flows as
\begin{subequations}\label{eq:coupling_conditions}
    \begin{align}
        f_\ell(t,L_\ell) = f_\ell^{out}(t)\\
        f_\ell(t,0) = f_\ell^{in}(t).\\
    \end{align}
\end{subequations}
Pipelines with compressors at their suction are denoted as active links $\ell \in \mathcal{L}_{a,gas}$.
The total power consumed by each compressor as part of each active link is then given by
\begin{equation}
    P_{\ell}(t)= c_p\cdot T\cdot f_\ell^{in}(t)\left(\left(\frac{p_\ell^{in}(t) + \Delta \Theta_{\ell}(t) }{p_\ell^{in}(t)}\right)^{\frac{\gamma-1}{\gamma}}-1\right), \quad \ell \in \mathcal{L}_{a,gas}, \label{eq:compression}\\
\end{equation}
where $c_p$, $T$, and $\gamma$ are constant parameters given in Table \ref{table:parameters}.
Pressure boundaries are defined by the following conditions:
\begin{align}
    p_\ell(t,L_\ell) = p_\ell^{out}(t) \quad \ell \in \mathcal{L}_{gas}\\
    p_\ell(t,0) = p_\ell^{in}(t) \quad \ell \in \mathcal{L}_{p,gas}\\
    p_\ell(t,0) = p_\ell^{in}(t) + \Delta \Theta_{\ell}(t) \quad \ell \in \mathcal{L}_{a,gas}.
\end{align}

Each pipeline also has a a steady-state initial condition:
\begin{subequations}\label{eq:weymouth_ss}
    \begin{align}
         c_{1,\ell}\frac{\partial f_\ell(0,x)}{\partial x} &= 0  \quad \ell \in \mathcal{L}_{gas}\label{eq:ssweymouthmass} \\
         c_{2,\ell}\frac{\partial p_\ell(0,x)}{\partial x} + c_{3,\ell}\frac{f_\ell(t,x)}{p_\ell(t,x)} \mods{f_\ell(t,x)} &= 0  \quad \ell \in \mathcal{L}_{gas},\label{eq:ssweymouthmomentum}
    \end{align}
\end{subequations}
and an operational constraint to return its linepack back to its starting inventory:
\begin{equation}
    m_{\ell}(T) \ge m_{\ell}(0)\label{eq:linepack_constraint} \quad \ell \in \mathcal{L}_{gas}.
\end{equation}
Here, the amount of linepack in any pipeline segment $m_\ell$ can be computed as:
\begin{equation}\label{eq:linepack}
    m_\ell(t) = \frac{1}{c_{1,\ell}}\int_{0}^{L_\ell}p_\ell(t,x)dx.
\end{equation}
Combining equations \eqref{eq:weymouth} through \eqref{eq:linepack_constraint} results in the following component model for each pipeline:
\begin{subequations}\label{eq:gas_network_model}
    \begin{align}
        &\frac{\partial p_\ell(t,x)}{\partial t} + c_{1,\ell}\frac{\partial f_\ell(t,x)}{\partial x} = 0 \quad \ell \in \mathcal{L}_{gas}\\
        &\frac{\partial f_\ell(t,x)}{\partial t} + c_{2,\ell}\frac{\partial p_\ell(t,x)}{\partial x} + c_{3,\ell}\frac{f_\ell(t,x)}{p_\ell(t,x)} \mods{f_\ell(t,x)} = 0 \quad ell \in \mathcal{L}_{gas}\\
        &p_{\ell}(L_{\ell},t)=\Theta_{rec(\ell)}(t),\quad \ell \in \mathcal{L}_{gas}\\
        &p_\ell(t,0) = p_\ell^{in}(t) \quad \ell \in \mathcal{L}_{p,gas}\\
        &p_\ell(t,0) = p_\ell^{in}(t) + \Delta \Theta_{\ell}(t) \quad \ell \in \mathcal{L}_{a,gas}\\
        &f_\ell(t,L_\ell) = f_\ell^{out}(t)\\
        &f_\ell(t,0) = f_\ell^{in}(t)\\
        &P_{\ell}(t)= c_p\cdot T\cdot f_\ell^{in}(t)\left(\left(\frac{p_\ell^{in}(t) + \Delta \Theta_{\ell}(t) }{p_\ell^{in}(t)}\right)^{\frac{\gamma-1}{\gamma}}-1\right), \quad \ell \in \mathcal{L}_{a,gas}\\
        &c_{1,\ell}\frac{\partial f_\ell(0,x)}{\partial x} = 0  \quad \ell \in \mathcal{L}_{gas} \\
        &c_{2,\ell}\frac{\partial p_\ell(0,x)}{\partial x} + c_{3,\ell}\frac{f_\ell(t,x)}{p_\ell(t,x)} \mods{f_\ell(t,x)} = 0  \quad \ell \in \mathcal{L}_{gas}\\
        &m_\ell(t) = \frac{1}{c_{1,\ell}}\int_{0}^{L_\ell}p_\ell(t,x)dx.
    \end{align}
\end{subequations}
These PDEs are discretized in space-time by using finite differences.

\subsection{Gas Junction Component  Model}

Each gas node model consists of a nodal pressure $\Theta_n$, a set of gas supplies $\mathcal{S}_n$, and a set of demands $\mathcal{D}_n$.
Gas nodes impose the following operational constraints on the network

\begin{subequations}\label{eq:node_model}
\begin{align}%\label{eq:gas_network_model}
    &\underline{\Theta}_n \le \Theta_n(t) \le \overline{\Theta}_n ,\quad n \in \mathcal{N}_{gas}\label{eq:node_pressure}\\
    &0 \le f_{n,d}(t) \le f_{n,d}^{target}(t)\label{eq:demand_target}, \quad d \in \mathcal{D}_n, \quad n \in \mathcal{N}_{gas}\\
    &0 \le f_{n,s}(t) \le \overline{s}_{n,s}(t)\label{eq:supply_limit}, \quad s \in \mathcal{S}_n, \quad n \in \mathcal{N}_{gas},
\end{align}
\end{subequations}

where $\underline{\Theta}_n$ is the lower pressure bound for the node, $\overline{\Theta}_n$ is the upper pressure bound,
$d_{n,d}(t)^{target}$ is the target demand for demand $d$ on node $n$ and $\overline{s}_{n,s}$ is the available gas generation
from supply $s$ on node $n$.

\subsection{Linking Constraints}
The linking constraints connect the pipeline and gas node models into a physical gas network
by enforcing boundary conditions and conservation constraints.
We can thus define the boundary conditions for the link pressures with
\begin{subequations}\label{eq:boundary_conditions}
    \begin{align}
        &p_{\ell}(L_{\ell},t)=\Theta_{rec(\ell)}(t), \quad \ell \in \mathcal{L}_{gas}\label{eq:bc1}\\
        &p_{\ell}(0,t)=\Theta_{snd(\ell)}(t), \quad \ell \in{\mathcal{L}_{gas}},\label{eq:bc2}
    \end{align}
\end{subequations}

and we can express the node balances around each gas node with the following link constraint
\begin{equation}\label{eq:node_balance}
    \sum_{\ell\in\mathcal{L}_n^{rec}}f^{out}_{\ell}(t)-\sum_{\ell \in\mathcal{L}_n^{snd}}f^{in}_{\ell}(t) +
    \sum_{s\in\mathcal{S}_n}f_{n,s}(t) - \sum_{d\in \mathcal{D}_n}f_{n,d}(t) = 0, \quad n \in\mathcal{N}_{gas}.
\end{equation}
We note that the supply flows $f_{n,s}(t) \quad s \in\mathcal{S}_n$ and the demand flows $f_{n,d}(t) \quad d \in \mathcal{D}_n$ are algebraic states of the system
(if the pressures of the corresponding nodes are specified) \cite{Zavala2014}.

\begin{table}[!h]
\caption{Variables and units for gas pipeline network model}
\begin{center}
\begin{tabular}{l l l}%{|l|l|}
\hline{Variable}       &  Description &  Units \\
\hline{$t$} & Time & $s$ \\
$x$ & Spatial dimension & $m$ \\
$\rho$ & Density & $kg / m^3$\\
$v$  & Velocity  & $m/s$\\
$V$  & Volumetric flow rate  & $m^3/s$\\
$p$  & Pressure  & $bar$\\
$f$  & Mass flow rate  & $kg/s$\\
$f^{in}$  & Pipe inlet flow  & $kg/s$\\
$f^{out}$  & Pipe outlet flow  & $kg/s$\\
$s$  & Supply flow  & $kg/s$\\
$d$  & Demand flow  & $kg/s$\\
$\Theta$  & Node Pressure  & $bar$\\
$\Delta \Theta$  & Boost Pressure  & $bar$\\
$P$  & Compressor Power  & $kW$\\
$m$  & Line-pack  & $kg$
\end{tabular}
\end{center}
\label{table:variables}
\end{table}

\begin{table}[!h]
\caption{Parameters and units for gas pipeline network model}
\begin{center}
\begin{tabular}{l l l}%{|l|l|}
\hline{Parameter}       &  Description &  Units \\
\hline{$\Delta x$,$\Delta t$} & Space and time discretization interval length & $m,s$ \\
$c$ & Speed of sound in gas & $m/s$ \\
$c_p,c_v$ &Gas heat capacity for constant pressure and volume & $2.34 kJ/kg K, 1.85 kJ/kg K$ \\
$\gamma, z$ & Isentropic expansion coefficient and compressibility & $-$ \\
$R$ & Universal gas constant & $8314 J / kmol K$ \\
$M$ & Gas molar mass & $18 kg/kmol$ \\
$\rho_n$ & Gas density at normal conditions& $0.72 kg/m^3$ \\
$T$ & Gas temperature & $K$ \\
$L,D,A$ & Pipeline length, diameter, and cross sectional area & $m,m,m^2$ \\
$\lambda,\epsilon$ & Friction factor and pipe rugosity & $-,m$ \\
$\alpha_f$ & Scaling factor for flow & $\frac{3600}{1 x 10^4 \rho_n} [=] \frac{SCMx10^4/h}{kg/s}$ \\
$\alpha_p$ & Scaling factor for pressure & $1x10^{-5} bar/ Pa$ \\
$c_1$ & Auxiliary constant & $\frac{bar/s}{SCMx10^4/h \ m}$ \\
$c_2$ & Auxiliary constant & $\frac{SCMx10^4/h - s}{bar/s}$ \\
$c_3$ & Auxiliary constant & $\frac{bar}{SCMx10^4 - s/h \ m}$ \\
$c_4$ & Auxiliary constant & $\frac{kW}{SCMx10^4/h}$
\end{tabular}
\end{center}
\label{table:parameters}
\end{table}

\FloatBarrier

\section{Benders Study Model}\label{app:benders}

The model for the two-stage stochastic resource management optimization problem is:
\begin{align}
    \min_{x,z,y} & \quad \sum_{\xi \in \Xi} \sum_{f \in \mathcal{F}} p_f(\xi) u_f(\xi) \label{eq:benders_problem}\\
    \textrm{s.t.} & \sum_{a \in \mathcal{A}_\mathcal{B}} c_a x_a + \sum_{j \in \mathcal{B}} h_j z_j \le Budget \nonumber\\
    &w_j = \gamma_j + z_j + \sum_{a \in rec(\mathcal{A}_{\mathcal{B}},j)} x_a - \sum_{a \in snd(\mathcal{A}_{\mathcal{B}},j)} x_a \quad j \in \mathcal{B} \nonumber\\
    &q_j(\xi) = w_j - \sum_{a \in snd(\mathcal{A}_{\mathcal{F}},j)} y_a(\xi) \quad j \in \mathcal{B} \nonumber\\
    &\sum_{a \in snd(\mathcal{A}_{\mathcal{F}},f)} y_a(\xi) + u_f(\xi) \ge d_f(\xi) \quad f \in \mathcal{F} \nonumber,
\end{align}

where $\Xi$ is a set of realized scenarios, $\mathcal{B}$ is a set of bases each containing resources, $\mathcal{F}$ is a set of districts with resource demands, $\mathcal{A}$ is a set
of arcs connecting bases and districts, $\mathcal{A}_{\mathcal{B}} \subseteq \mathcal{A}$ is the set of arcs between bases, and $\mathcal{A}_{\mathcal{F}}$ is the
set of arcs between bases and districts. Parameter $p_f(\xi)$ is the probability that scenario $\xi$ realizes at district $f$ and variable $u_f(\xi)$ is the
unmet demand at district $f$ after dispatch decisions are made for scenario $\xi$. Variable
$x \in \mathbb{R}^{|\mathcal{A}_\mathcal{B}|}$ is a first stage decision to move resources between bases, $z \in \mathbb{R}^{|\mathcal{B}|}$ is a
first stage decision to purchase resources at bases, and variable $y \in \mathbb{R}^{|\mathcal{A}_\mathcal{F}|}$ is a second stage decision
to dispatch resources to districts after realizing district demands.  Parameter $\gamma_j$ is the initial amount of resources in base $j$, $w_j$ is the amount of resources
in base $j$ after making transfers, variable $q_j(\xi)$ is the amount of resources in each base after dispatching to districts for scenario $\xi$, and parameter $d_f(\xi)$ is the
resource demand of district $f$ for scenario $\xi$.

Problem \eqref{eq:benders_problem} is reformulated to conduct Benders decomposition by decomposing it into a master problem and subproblem for each scenario.  The master problem is given by:
\begin{align}\label{eq:benders_master}
    \min_{x,z,\theta_{cut}} & \ \qquad \theta_{cut}\\
    s.t. & \sum_{a \in \mathcal{A}_\mathcal{B}} c_a x_a + \sum_{j \in \mathcal{B}} h_j z_j \le Budget \nonumber\\
    &w_j = \gamma_j + z_j + \sum_{a \in rec(\mathcal{A}_{\mathcal{B}},j)} x_a - \sum_{a \in snd(\mathcal{A}_{\mathcal{B}},j)} x_a, \quad j \in \mathcal{B} \nonumber\\
    &x_a \ge 0, \quad a \in \mathcal{A} \nonumber \\
    &z_j \ge 0, \quad j \in \mathcal{B} \nonumber \\
    &w_j \ge 0, \quad j \in \mathcal{B} \nonumber \\
    &\theta_{cut} \ge 0 \nonumber\\
    &\theta_{cut} \ge c(x), \quad c \in C, \nonumber
\end{align}
where $\theta_{cut}$ is a variable to enforce cutting planes and $C$ is the set of cutting planes added to the master problem.
The subproblem is a function of the master solution variable $\hat{w}$ and scenario $\xi$ and is given by \eqref{eq:benders_sub}.
\begin{align}\label{eq:benders_sub}
     Q(\hat{w},\xi) &:= \min_{y} \quad \sum_{f \in \mathcal{F}} p_f(\xi) u_f \\
     \textrm{s.t.} \quad & q_{j} =  \hat{w}_j - \sum_{a \in snd(\mathcal{A}_{\mathcal{F}},j)} y_a, \quad j \in \mathcal{B} \nonumber\\
     & \sum_{a \in snd(\mathcal{A}_{\mathcal{F}},j)}y_a + u_f \ge d_f(\xi), \quad f \in \mathcal{F} \nonumber\\
     &q_j \ge 0, \quad j \in \mathcal{B} \nonumber\\
     &y_a \ge 0, \quad a \in \mathcal{A}_{\mathcal{F}} \nonumber\\
     &u_f \ge 0, \quad f \in \mathcal{F}. \nonumber
\end{align}
The master node $m$ in the computing graph contains two tasks described by Task \ref{task:run_master}
and Task \ref{task:receive_solution}.  Each subnode in the computing graph consists of a single task described by Task \ref{task:run_sub_problem}.
Task \ref{task:run_master} runs the master problem \eqref{eq:benders_master} and updates the first stage solution $\hat{x}$ (which contains all first-stage variables including $\hat{w}$) and checks whether
convergence has been satisfied.  If not satisfied, it updates the master node scenario attributes $\{\xi_1,...,\xi_N\}$ with the first $N$ values from the
scenario set $\Xi$.  Task \ref{task:receive_solution} runs when the master node receives an update to a solution attribute $s_n$.  The task checks
whether every solution has returned, and if true, it updates the flag attribute $flag$, which indicates that the master problem is ready to be solved.  If not all
subproblem solutions have returned, the task updates the attribute $\xi_n$ with the next scenario in $\Xi$.  Task \ref{task:run_sub_problem} solves the subproblem \eqref{eq:benders_sub} given
a first stage solution $\hat{x}$ and a scenario $\xi$ and updates the subnode solution attribute $s$.

\makeatletter
\renewcommand*{\ALG@name}{Task}
\makeatother
\setcounter{algorithm}{0}
\begin{algorithm*}[h!]
\caption{${\tt run\_master}$\;(computing graph $\mathcal{G}$,node $m$)}
\label{task:run_master}
\begin{algorithmic}[1]
    \State Get node attribute $C$: Set of current master problem cuts
    \State Get node attribute $S$: Set of solutions received from sub-problems
    \State Get node attribute $\Xi$: Set of sample scenarios
    \State Get node attribute $\hat{x}$: Master problem solution attribute
    \State Get node attributes $\{\xi_{1},...,\xi_{N}\}$: Scenarios to subnodes $\mathcal{N}$
    \State Solve \eqref{eq:benders_master}: ${\tt \bm{update \ attribute}} \ \hat{x} \gets {\tt solve\_master\_problem}(C)$
    \If{${\tt problem\_converged}$($\hat{x},S$):}
    \State Stop Computing Graph $\mathcal{G}$
    \Else
    \State Start sending new scenarios: ${\tt \bm{update \ attributes}} \{\xi_{1},...,\xi_{N}\} \gets \Xi[1:N]$
    \EndIf
    \State Empty current solutions: ${\tt set} \ S \gets \{\}$
\end{algorithmic}
\end{algorithm*}\label{alg:master_node_tasks}

\begin{algorithm*}[h!]
\caption{${\tt receive\_solution}$(node $m$, attribute $s_{n}$)}\label{task:receive_solution}
\begin{algorithmic}[1]
    \State Get node attribute $S$: Set of solutions received from sub-problems
    \State Get node attribute $\Xi$: Set of sample scenarios
    \State Get argument $s_{n}$: Solution from subnode $n$
    \State Get node attribute $\xi_{n}, \; n \in \mathcal{N}$: Subnode $n$ scenario attribute
    \State Get node attribute $flag$: Flag that master problem is ready to solve
    \State ${\tt push \ to} \ S \gets {\tt get\_objective\_value}(s_{n})$
    \State ${\tt push \ to} \ C\gets {\tt compute\_new\_cut}(s_{n})$
    \If{$\tt{all\_scenarios\_complete}$($S$)}
    \State \text{Master problem is ready}: ${\tt \bm{update \ attribute}} \ flag$
    \Else
    \State \text{Send new scenario to node} $n$ : {${\tt \bm{update \ attribute}} \ \xi_{n} \gets {\tt new\_scenario}(\Xi)$}
    \EndIf
\end{algorithmic}
\end{algorithm*}\label{alg:master_node_tasks}

\begin{algorithm*}[h!]
\caption{{\tt solve\_subproblem}(node $n$)}
\label{task:run_sub_problem}
\begin{algorithmic}[1]
    \State Get node attribute $\hat{x}$: Master solution received on subnode $n$
    \State Get node attribute $\xi$: Scenario received on subnode $n$
    \State Solve \eqref{eq:benders_sub}: ${\tt update \ attribute} \ s \gets {\tt solve\_subproblem}(\hat{x},\xi)$
\end{algorithmic}
\end{algorithm*}\label{alg:master_node_tasks}

\section{MPC Study Model}\label{app:reactor}

The model for the plant is given by the following set of differential equations:
\begin{align}\label{eq:plant}
    &\frac{dH_1}{dt} = \frac{1}{\rho A_1}(F_{f1} + F_R - F_1)\\
    &\frac{dx_{A1}}{dt} = \frac{1}{\rho A_1 H_1}(F_{f1} x_{A0} + F_R x_{AR} - F_1 x_{A1}) - k_{A1} x_{A1} \nonumber \\
    &\frac{dx_{B1}}{dt} = \frac{1}{\rho A_1 H_1}(F_R x_{BR} - F_1 x_{B1}) + k_{A1} x_{A1} - k_{B1} x_{B1} \nonumber \\
    &\frac{dT_1}{dt} = \frac{1}{\rho A_1 H_1}(F_{f1} T_0 + F_R T_R - F_1 T_1) - \frac{1}{C_p}(k_{A1} x_{A1} \Delta H_A + k_{B1} x_{B1} \Delta H_B) + \frac{Q_1}{\rho A_1 C_p H_1} \nonumber \\
    &\frac{dH_2}{dt} = \frac{1}{\rho A_2}(F_{f2} + F_1 - F_2)\nonumber \\
    &\frac{dx_{A2}}{dt} = \frac{1}{\rho A_2 H_2}(F_{f2} x_{A0} + F_1 x_{A1} - F_2 x_{A2}) - k_{A2} x_{A2} \nonumber \\
    &\frac{dx_{B2}}{dt} = \frac{1}{\rho A_2 H_2}(F_1 x_{B1} - F_2 x_{B2}) + k_{A2} x_{A2} - k_{B2} x_{B2} \nonumber \\
    &\frac{dT_2}{dt} = \frac{1}{\rho A_2 H_2}(F_{f2} T_0 + F_1 T_1 - F_2 T_2) - \frac{1}{C_p}(k_{A2} x_{A2} \Delta H_A + k_{B2} x_{B2} \Delta H_B) + \frac{Q_2}{\rho A_2 C_p H_2} \nonumber \\
    &\frac{dH_3}{dt} = \frac{1}{\rho A_3}(F_2 - F_D - F_R - F_3)\nonumber \\
    &\frac{dx_{A3}}{dt} = \frac{1}{\rho A_3 H_3}(F_2 x_{A2} - (F_D + F_R) x_{AR} - F_3 x_{A3}) \nonumber \\
    &\frac{dx_{B3}}{dt} = \frac{1}{\rho A_3 H_3}(F_2 x_{B2} - (F_D + F_R) x_{BR}) - F_3 x_{B3}) \nonumber \\
    &\frac{dT_3}{dt} = \frac{1}{\rho A_3 H_3}(F_2 T_2 - (F_D + F_R) T_R - F_3 T_3) + \frac{Q_3}{\rho A_3 C_p H_3}, \nonumber
\end{align}
where for $i = \{1,2\}$ we have:
\begin{equation*}
    k_{Ai} = k_A \exp \Big(-\frac{E_A}{R T_i} \Big), \quad k_{Bi} = k_B \exp \Big(-\frac{E_B}{R T_i} \Big).
\end{equation*}
The recycle and weight fractions are given by:
\begin{align*}
    &F_D = 0.01 F_R, \quad x_{AR} = \frac{\alpha_A x_{A3}}{\overline{x}_3}, \quad x_{BR} = \frac{\alpha_B x_{B3}}{\overline{x}_3}\\
    &\overline{x}_3 = \alpha_A x_{A3} + \alpha_B x_{B3} + \alpha_C x_{C3}, \quad x_{C3} = (1 - x_{A3} - x_{B3}).
\end{align*}
The target steady-state (set-point) is described by the parameters in Table \ref{table:reactor_parameters} and the initial
operating condition is given in Table \ref{table:reactor_initial}.
\begin{table}[h]
\caption{Target steady-state and parameters for reactor-separator system}
\begin{center}
\begin{tabular}{l l l l l l}%{|l|l|}
\hline
{Parameter}  &  Value &  Units & Parameter & Value & Units \\
\hline
$H_1$    & 16.1475 & $m$     & $A_1$    & 1.0   & $m^2$\\
$x_{A1}$ & 0.6291  & $wt$    & $A_2$    & 1.0   & $m^2$\\
$x_{B1}$ & 0.3593 & $wt$    & $A_3$    & 0.5   & $m_2$\\
$T_1$    & 387.594   & $K$     & $\rho$   & 1000   & $kg/m^3$\\
$H_2$    & 12.3137   & $m$     & $C_p$    & 4.2   & $kJ/kg K$\\
$x_{A2}$ & 0.6102 & $wt$    & $x_{A0}$ & 0.98 & $wt$ \\
$x_{B2}$ & 0.3760 & $wt$    & $T_0$    & 359.1   & $K$\\
$T_2$    & 386.993  & $K$     & $k_A$    & 2769.44  & $1/s$\\
$H_3$    & 15.0   & $m$     & $k_B$    & 2500.0   & $1/s$\\
$x_{A3}$ & 0.2928 & $wt$    & $E_A / R$    & 6013.95   & $kJ/kg$\\
$x_{B3}$ & 0.67 & $wt$    & $E_B / R$    & 7216.74   & $kJ/kg$\\
$T_3$    & 387.01 & $K$     & $\Delta H_A$ &  -167.4 & $kJ/kg$\\
$F_{f1}$ & 6.3778 & $kg/s$  & $\Delta H_B$ & -139.5 & $kJ/kg$\\
$Q_1$    & 26.0601 & $kJ/s$  & $\alpha_A$   & 5.0 &  --\\
$F_{f2}$ & 6.8126 & $kg/s$  & $\alpha_B$   & 1.0 &  --\\
$Q_2$    & 5.0382 & $kJ/s$  & $\alpha_C$   & 0.5 &  --\\
$F_R$    & 56.7989 & $kg/s$  &              &   &\\
$Q_3$    & 5.0347 & $kJ/s$  &              &   &\\
$F_1$    & 63.1766 & $kg/s$  &              &   &\\
$F_2$    & 69.9892 & $kg/s$  &              &   &\\
$F_3$    & 12.6224 & $kg/s$  &              &   &
\end{tabular}
\end{center}
\label{table:reactor_parameters}
\end{table}

\begin{table}[h]
\caption{Initial conditions for simulation of reactor-separator system}
\begin{center}
\begin{tabular}{l l l l l l}%{|l|l|}
\hline
State  &  Value &  Units & Input & Value & Units \\
\hline
$H_1$    & 25.4702 & $m$     & $F_{f1}$   & 1.1866    & $kg/s$\\
$x_{A1}$ & 0.1428 & $wt$    & $Q_1$    & 29.0597   & $kJ/s$\\
$x_{B1}$ & 0.7045 & $wt$    & $F_{f2}$    & 7.0263   & $kg/s$\\
$T_1$    & 415.944  & $K$     & $Q_2$   & 5.1067   & $kJ/s$\\
$H_2$    & 5.4703  & $m$     & $F_R$    & 11.6962    & $kg/s$\\
$x_{A2}$ & 0.3653 & $wt$    & $Q_3$ & 5.09834  &  $kJ/s$ \\
$x_{B2}$ & 0.5307 & $wt$    &  $F_1$  & 12.8828  & $kg/s$\\
$T_2$    & 399.303 & $K$     &  $F_2$  &  19.9091  & $kg/s$\\
$H_3$    & 15.0  & $m$     &  $F_3$  & 8.0960  & $kg/s$\\
$x_{A3}$ & 0.1565 & $wt$     \\
$x_{B3}$ & 0.67 & $wt$     \\
$T_3$    & 399.364 & $K$    \\
\end{tabular}
\end{center}
\label{table:reactor_initial}
\end{table}
\FloatBarrier

\noindent For the 3 MPC controllers, we have the following outputs and inputs:

\begin{align*}
    &y_1 = [H_1 ,\ x_{A1}, \ x_{B1}, \ T_1] ,\quad u_1 = [F_{f1}, \ Q_1, \ F_1]\\
    &y_2 = [H_2, \ x_{A2}, \ x_{B2}, \ T_2],\quad u_2 = [F_{f2}, \ Q_2, \ F_2]\\
    &y_3 = [H_3, \ x_{A3}, \ x_{B3}, \ T_3] ,\quad  u_3 = [F_R, \ Q_3, \ F_3].
\end{align*}
Each MPC controller uses a quadratic cost function with weights:
\begin{align*}
    &Q_{y1} = \textrm{diag}(100,10,100,0.1) \quad Q_{y2} = \textrm{diag}(10,10,100,0.1) \quad Q_{y3} = \textrm{diag}(1,10,10^5,0.1)\\
    &R_{yi} = \textrm{diag}(100,100,100),\;  i = \{1,2,3\}.
\end{align*}
The differential equations are discretized using an Euler scheme with a time horizon of $N = 20$ and a time step $\Delta \ t = 30$.
The cooperative MPC computation tasks are defined in Tasks \ref{task:run_plant}, \ref{task:control_action}, and \ref{task:receive_policy}.
Task \ref{task:run_plant} simulates the plant forward in time from the current computing graph time to the time of the next signal in the computing graph queue and updates the attribute $x$.
Task \ref{task:control_action} computes the open-loop control trajectory for MPC controller $n$ and updates its control injection $u_{inject}$ if it has completed enough iterations and updates
its control policy $u_{p_n}$ if it has not.  Task \ref{task:receive_policy} checks whether the MPC controller has recieved updates from the other MPC controllers and updates the flag indicator $flag$
if it has received both inputs.

\begin{algorithm*}[h!]
\caption{${\tt run\_plant}$(computing graph $\mathcal{G}$,node $n$)}
\label{task:run_plant}
\begin{algorithmic}[1]
    \State Get graph clock value $t_{now}$: Current clock time
    \State Get graph clock value $t_{next}$: Next signal time
    \State Get node attributes $\{u_1,u_2,u_3\}$: Received MPC controller actions
    \State Get node attribute $x$: Plant state
    \State Simulate plant \eqref{eq:plant} ${\tt update \ attribute} \ x \gets$ ${\tt simulate\_plant}$($u_1,u_2,u_3,x,t_{now},t_{next}$)
\end{algorithmic}
\end{algorithm*}

\begin{algorithm*}[h!]
\caption{${\tt control\_action}$(node $n$)}
\label{task:control_action}
\begin{algorithmic}[1]
    \State Get node attribute $iter_{max}$
    \State Get node attribute $y$: Received plant measurement
    \State Get node attribute $u_{inject}$: Current injected control from MPC controller
    \State Get node attributes $\{u_{p_i}, i \ne n\}$: Received control actions from other MPC controllers
    \State Get node attribute $iter$: Current cooperative iteration
    \State Calculate control action: $u_{calc} \gets {\tt compute\_control\_action}(y_n,u_{p_i}, i \ne n)$
    \State Update iteration count: $iter \gets iter + 1$
    \If{$iter = iter_{max}$}
    \State Update injected control:  ${\tt update \ attribute} \ u_{inject} \gets u_{calc}$
    \State $iter \gets 1$
    \Else
    \State Update control: ${\tt update \ attribute} \ u_{p_n} \gets u_{calc}$
    \EndIf
\end{algorithmic}
\end{algorithm*}

\begin{algorithm*}[h!]
\caption{${\tt receive\_policy}$(node $n$, attribute $u_{p_i}$)}
\label{task:receive_policy}
\begin{algorithmic}[1]
    \State Get node argument $u_{p_i}$ : Received neighbor controller actions
    \State Get node attribute $flag$ : Flag that control action is ready to compute
    \If{Received both policies( $ \{ u_{p_i}, i \ne n \} $ )}
    \State ${\tt update \ attribute} \ flag$
    \EndIf
\end{algorithmic}
\end{algorithm*}

\bibliography{plasmo}

\end{document}